\documentclass[11pt]{article}
\usepackage{amsfonts}
\usepackage{amsmath,amsthm,amscd,amssymb,mathrsfs,setspace, textcomp}

\usepackage{latexsym,epsf,epsfig,float}
\usepackage{color}
\usepackage[hmargin=2.25cm,vmargin=2.75cm]{geometry}

\usepackage{hyperref}
\usepackage{backref}

\newcommand{\ds}{\displaystyle}

\newcommand{\cD}{\mathscr{D}}

\theoremstyle{plain}
\newtheorem{theorem}{Theorem}[section]

\newtheorem{definition}{Definition}

\theoremstyle{remark}
\newtheorem{remark}{Remark}[section]
\numberwithin{equation}{section}
\numberwithin{theorem}{section}
\numberwithin{remark}{section}
\numberwithin{assumption}{section}
\numberwithin{condition}{section}

\title{A Thorough Look at the (In)Stability of Piston-theoretic Beams}

 \author{\normalsize \begin{tabular}[t]{c@{\extracolsep{.6em}}c@{\extracolsep{.6em}}c@{\extracolsep{.6em}}c}
      Jason Howell\footnote{Carnegie Mellon University, PA}&  ~~Katelynn Huneycutt\footnote{University of Maryland, Baltimore County, MD}& ~~Justin T. Webster\footnote{University of Maryland Baltimore County, MD}&~~ Spencer Wilder\footnote{College of Charleston, SC}\\
\it howelljs@cmu.edu & \it dd31975@umbc.edu &  \it websterj@umbc.edu & \it wildersb@g.cofc.edu 
\end{tabular}}

\begin{document}
\maketitle

\begin{abstract} {\noindent We consider a beam model representing the transverse deflections of a one dimensional elastic structure immersed in an axial fluid flow. The model includes a nonlinear elastic restoring force, with damping and non-conservative terms provided through the flow effects. Three different configurations are considered: a clamped panel, a hinged panel, and a flag (a cantilever clamped at the leading edge, free at the trailing edge). After providing the functional framework for the dynamics, recent results on well-posedness and long-time behavior of the associated dynamical system for solutions are presented. Having provided this theoretical context, in-depth numerical stability analyses are provided, focusing both at the onset of flow-induced instability (flutter), and qualitative properties of the post-flutter dynamics {\em across configurations}. Modal approximations are utilized, as well as finite difference schemes.  
  \\[.3cm]
\noindent {\bf Key terms}:  extensible beam, nonlinear elasticity, flutter, stability, attractors, modal analysis
 \\[.3cm]
}
\end{abstract}

\maketitle

\section{Introduction}

In this treatment we provide a multifaceted numerical study---{from the dynamical systems point of view}---for a simplified class of 1-D models representing the physical phenomenon known as {\em aeroelastic flutter}. Flutter is fluid-structure (feedback) instability that occurs between elastic displacements of a structure and responsive aerodynamic pressure changes at the flow-structure interface \cite{bolotin,dowell,survey2}. The onset of flutter, for particular flow parameters, represents a bifurcation in the linearization of the system \cite{holmes}, and the qualitative properties of the post-flutter dynamics can be analyzed from an infinite dimensional/control-theoretic point of view \cite{HolMar78,survey2,fereisl}.  

Beam (or plate) flutter is a topic of great interest in the engineering literature \cite{bolotin,dowellnon,pist2,b:paidoussis:98,inext2,paidoussis3,jfs,vedeneev}, as well as, more recently, the mathematical literature \cite{survey2, HLW,conequil1,conequil2} (and many references therein). Though there is a vast (mostly engineering-oriented) literature on flow-structure instabilities, we provide a straight-forward analysis utilizing a 1-D structural model in multiple configurations of interest. In particular, we consider both linear and nonlinear extensible beam models with piston-theoretic flow terms that, although simplified, demonstrate good correspondence of the dynamics with what is empirically known about the onset of flutter and post-flutter behaviors \cite{jfs,pist2}. We will demonstrate the remarkable ability of these simplified equations to capture a robust set of un/stable behaviors associated with structures in axial flow. We also provide modern, rigorous context for these numerical studies by discussing recent theoretical results for the nonlinear dynamical systems associated to the PDE models given here.

Our studies primarily take the form of computational stability analyses for a beam distributed on $x \in [0,L]$, immersed in an inviscid potential flow of sufficiently high Mach number \cite{pist2}; we consider three beam configurations: clamped-clamped panel ({\bf C}), hinged-hinged panel ({\bf H}), and clamped-free {\em cantilever} ({\bf CF}). The fluid flows with steady velocity $U \ge 0$ in the $x$-direction (along the length of the beam), in contrast with so-called normal flow; thus, $x=0$ is the {\em leading edge}, and $x=L$ is the {\em trailing edge}. The cantilever {\bf CF} is taken to be a {\em flag}, clamped at the leading edge, and free at the  trailing edge\footnote{A recent configuration of interest is the so called {\em inverted flag}, where the free edge leads; see \cite{amjad, inverted1}}. The baseline linear structural model we consider is the standard Euler-Bernoulli beam, though for the {\bf CF} configuration we mention the so called Rayleigh beam (allowing for rotational inertia).\footnote{See \cite{slenderness} for a nice description and comparison of four fundamental beam theories.} Additionally, we consider {\em extensible nonlinear beams} in sense of Krieger-Woinowsky \cite{wonkrieg}---see \cite{HTW} for more details. For the flow theory, we employ the {\em linear piston theory}, which reduces the flow dynamics to a simple non-conservative right hand side written through the material derivative of structural variable \cite{HolMar78,pist2}.

As with all {flutter problems} \cite{dowell,b:paidoussis:98}, and as we will demonstrate in the sequel, the {\em onset} of instability can be studied via linear theory---typically as an eigenvalue problem (see e.g., \cite{jfs,vedeneev}). Indeed, if one can model the dynamic load across the structure, the issue is to determine under what conditions the non-conservative flow loading destabilizes structural (eigen)modes. However, to study the qualitative properties of the dynamics in the post-flutter regime, analysis will require {\em some} physical nonlinear restoring force to keep trajectories {bounded in time} \cite{dowellnon, fereisl,survey2}. Indeed, for the post-flutter dynamics of a given system, the theory of nonlinear dynamical systems is often used, e.g.,  compact global {\em attractors} \cite{springer, survey2} as well as {\em exponential attractors} \cite{HLW,quasi,springer}.

\subsection{Focus of This Study and Relation to Previous Literature} 

This study considers both pre-~and post-flutter extensible, piston-theoretic beams in three configurations, from two numerical points of view: (i) a ``modal" analysis, and (ii) a traditional finite difference scheme in space, both described in detail below.

We consider the effect of the piston-theoretic flow terms on the linear and nonlinear beam. We offer a thorough study of each key parameter's effect on stability, in the sense of: convergence to an equilibrium, onset of dynamic instability, convergence to a limit cycle oscillation (or chaos), and qualitative properties of non-transient dynamics. The analysis spans the three most physically relevant configurations appearing in the literature, and, unlike most previous work (e.g., \cite{holmes,vedeneev,modalcant}), we provide comparisons of the dynamics across the configurations. The inclusion of structural nonlinearity distinguishes this work from \cite{vedeneev,jfs,shubov1,modalcant}, allowing us to discuss qualitative properties of post-flutter behaviors. Our detailed parameter studies are done while bearing in mind the large body of {\em recent theoretical results}, in particular those in \cite{HLW,survey2}, and those in \cite{HTW} for cantilevers with piston-theoretic terms. 

Our analysis is done in the spirit of the impressive (by now classical) references \cite{holmes,HolMar78} that study a piston-theoretic 1-D, {\bf H} panel, allowing for in-plane pre-stressing. These papers consider various unstable behaviors, including divergence/buckling (static instability) and flutter (dynamic instability): the earlier \cite{holmes} performs a tremendously thorough ODE analysis for a truncated---low mode---system, while \cite{HolMar78} provides an infinite dimensional framework for discussions instability. We also mention the classical reference \cite{dowellchaos}, where the same pre-stressed panel model is reduced to a two degree of freedom dynamical system via spectral methods, and considered as a deterministic example exhibiting chaos (for sufficiently large flow and stressing parameters). 
 
For more general studies of fluttering beams and plates (with extensive references), we point to \cite{dowell}.  In the discussion of piston theory, we mention the modern engineering references \cite{jfs,vedeneev}, classic engineering references \cite{pist2,bolotin}, and the recent mathematical-minded survey \cite{survey2}. Other mathematical references addressing piston-theoretic models include \cite{HLW, Memoires, daniellorena,springer}. The article \cite{HLW}: (i) performs numerical simulations (which guide the simulations herein); (ii) theoretically addresses various classes of {\em attractors} arising in the {\bf C} beam/plate configuration.  We follow much of that analysis there, in particular tracking stability/instability of dynamics with respect to piston-theoretic terms, as well as types and size of damping effects. The recent \cite{HTW} follows the general theory/numerics approach in \cite{HLW}, but focuses specifically on the {\bf CF} configuration and the effect of rotational inertia in the beam.

\subsection{Modeling and Discussion}\label{modeling}
We now address the modeling  pertinent to the theoretical discussion and numerical studies below.

Though there are various ways to consider flow-beam coupling, the simplest is to eliminate the fluid dynamic variables altogether. Such a tack has the benefit of reducing stability questions to a single, non-conservative dynamics. Additionally, the reduction eliminates moving boundary issues for the flow (particularly for the cantilever). On the other hand, such a reduction is a dramatic simplification of complex, multi-physics phenomena;  however, focusing on the simple model allows us use robust mathematical theory and perform a thorough numerical study that can be exposited straight-forwardly. (More sophisticated flow-structure models are certainly explored in the rigorous mathematical literature \cite{delay,webster,springer,fereisl}).

Motivated by the bulk of engineering studies cited thus far, we consider beam dynamics interacting with a potential flow. For high Mach numbers, the dynamic pressure on the surface of the structure, given by $p(x,t)$, can be approximated point-wise in $x$ by the pressure on the head of a piston moving through a column of fluid \cite{vedeneev,pist2}. The pressure can then be written in terms of the {\em down-wash} of the fluid $W=(\partial_t+U\partial_x)w$, where $w(x,t)$ is the transverse displacement of the beam, and $U$ is the unperturbed axial flow velocity. This results in a nonlinear expression \cite{pist2} in  $W$ $$\ds   p_*(x,t) = p_0 - \mu \left[1+\dfrac{w_t+Uw_x}{U} \right]^{\gamma},~~\gamma \approx 1.4~~ \cite{dowell},$$  which is then linearized to produce the so called {\em linear piston theory} used on the RHS of the beam:
\begin{equation}\label{linpist}
p(x,t)=p_0(x)-\beta[w_t+Uw_x], ~~\text{for}~U>\sqrt2.
\end{equation}
Above, $p_0(x)$ is a static pressure on the surface of the beam (for the numerical portion of this paper we will take $p_0(x)=0$). The parameter $\beta>0$ is a fluid density parameter that typically depends on $U$, though numerically we consider $\beta$ decoupled from $U$.\footnote{See \cite{vedeneev} for a discussion of the flow non-dimensionalization, and further discussion of characteristic parameter values.}

For the structural model, we assume a slightly extensible beam \cite{lagleug}. Taking the standard Kirchhoff-type hypotheses \cite{slenderness,lagnese} leads to an Euler-Bernoulli beam, or the so called Rayleigh beam when allowing for rotational inertia in beam filaments \cite{slenderness}. For nonlinear effects, we invoke a {\em large deflection} elasticity model, taking into account quadratic terms in the strain-displacement relation \cite{lagnese}. From this, we obtain the beam analog of the full von Karman plate equations \cite{lagleug}, that take into account both in-plane $u$ and out-of-plane $w$ (Lagrangian) dynamics. The  principal nonlinear effect accounts for {\em the local effect of stretching on bending} (thus requiring beam extensibility). 
Let $(w,w_t)$ correspond to the out-of-axis (transverse) dynamics, and $(u,u_t)$ the in-axis (longitudinal) dynamics. Rotational inertia in  beam filaments is represented by $\alpha \ge 0$, and $D_1,D_2 \ge 0$ are elastic coefficients.\footnote{The coefficients $\alpha$, $D_1$, $D_2$ are related as $\alpha D_1= D_2$ in the physical derivation of the beam equation---see \cite{lagleug}.} 
 \begin{equation}\label{LLsystem}
\begin{cases} u_{tt} -D_1\left[u_x+\frac{1}{2}(w_x)^2\right]_x=0 \\
(1-\alpha\partial_x^2)w_{tt}+D_2\partial_x^4 w -D_1\left[ w_x(u_x+\frac{1}{2}w_x^2)\right]_x = p(x,t) \\
u(t=0)=u_0;~~u_t(t=0)=u_1;~~w(t=0)=w_0;~~w_t(t=0)=w_1 \\
BC(u,w). \end{cases}
\end{equation}
The boundary conditions $BC(u,w)$ are determined by the configuration to be studied.
\begin{remark}
At present, all rigorous mathematical analyses of the conservative \eqref{LLsystem} model (or the 2-D analog) {\em require} $\alpha>0$ \cite{lagleug,HTW}. However, when considering purely transverse $w$ dynamics of thin beams and plates, $\alpha=0$ is almost always taken (see the discussion in \cite[Ch.1]{lagnese}).\end{remark}
The beam model under consideration here is a simplification of the above system \eqref{LLsystem}:  take $u_{tt} \approx 0$ and solve the resulting equation for $u$ in terms of $w$, reducing the system to transverse dynamics.
In eliminating the $u$ dynamics, we enforce the condition that $u(L)$ is zero---an unproblematic assumption for {\bf C} and {\bf H}, but certainly an imposition for {\bf CF} (see Remark \ref{assumption}). A non-zero choice corresponds (mathematically) to pre-stressing the beam at equilibrium. 
Completing the simplification, one obtains the beam dynamics considered here:
\begin{equation}\label{Bergerplate*}
(1-\alpha\partial_x^2)w_{tt}+D\partial_x^4 w +[b_1-b_2||w_x||^2]w_{xx} = p_0(x)-\beta[w_t+Uw_x] \\
\end{equation}
Above, we renamed or introduced various parameters: $b_1 \in \mathbb R$ measures in-plane stretching ($b_1<0$) or compression ($b_1>0$) at equilibrium (pre-stressing); $b_2>0$ gives the strength of the effect of stretching on bending (i.e., the nonlinear restoring force). We forgo a detailed discussion of the inertial parameter $\alpha \ge 0$, as well as associated structural damping in the model, until the next section.
\begin{remark}
This model has been historically referred to as the Krieger-Woinowsky beam \cite{wonkrieg, beam4}, but has also been referred to as Kirchhoff or Berger---the latter used owing to the plate equation of the same nonlinear structure \cite{HLW,Memoires}.  Classic references \cite{ball, ball2,biacrip,dickey, edenmil} discuss well-posedness and long-time behavior for extensible beams and Berger plates, typically with {\bf C} or {\bf H}  boundary conditions.
\end{remark}
The notation $BC(w)$ will provide the boundary conditions corresponding to the pertinent configuration  ({\bf C}, {\bf H}, or {\bf CF}) given mathematically in the next section. These are standard boundary conditions, {\em except in the} {\bf CF} {\em case} when imposing $u(L)\approx 0$ necessitates a new condition: 
\begin{align}\label{badBC}
~-\alpha\partial_xw_{tt}(L)+Dw_{xxx}(L)+[b_1-b_2||w_x||^2]w_x(L)=0.
\end{align}
It is certainly noteworthy that the inherited boundary condition is nonlinear and nonlocal (but collapses to the standard free boundary condition if $b_1=b_2=0$). \eqref{badBC} is mathematically non-trivial, as well as numerically challenging---\cite{HTW} addresses this condition, and we elaborate on it below in Section \ref{postflutter}. 
\begin{remark}\label{assumption} For the cantilever, the restriction $u(L)=const.$ throughout deflection is not entirely physical, at least in the absence of external mechanical effects. Indeed, \eqref{badBC} implicitly confines the free end of the beam to a fixed vertical line throughout deflection, thereby {\em inducing} stretching that provides the nonlinear restoring force. For a post-flutter cantilever, free end displacements can be quite large, with non-trivial inertial effects \cite{inext1,inext2}. However, in using piston theory we assume larger values of $U$, known to reduce the magnitude of $u(L)$, thereby improving the nonlinear structural model's viability.

With this in mind, a ``good"  cantilever flutter model should allow for in-plane displacements, and should not  primarily consider extensible effects. However, the added complexity of the nonlinearly coupled $u$/$w$ system in  \eqref{LLsystem} or the recent {\em inextensible beam} models \cite{inext1,inext2,paidoussis,paidoussis3} dramatically increases computational and analytical difficulties. Rigorous results for such models are few.  
\end{remark}
\noindent {\bf In summary}: Our studies here provide a thorough analysis of simple, nonlinear fluttering beams, across various configurations. We acknowledge that the model considered here is greatly simplified: (i) it is only 1-D, not accounting for in-axis dynamics; (ii) nonlinear effects studied here are strictly due to extensibility, and we do not consider more sophisticated nonlinear models; (iii) we do not fully model flow dynamics, instead relying on piston-theoretic simplifications. Yet, even this simplified model yields interesting behaviors worthy of rigorous study, complemented by a thorough numerical examination done from a dynamical systems point of view. Our conclusions verify empirical conclusions from the engineering literature, as well as provide validation for recent abstract (infinite dimensional) results for beams and plates \cite{HLW,HTW}. Our studies also serve as a  baseline for future studies of more complex models, testing conjectured similarities between these dynamics in certain regimes.

\section{PDE Model and Configurations} 

Recalling from Section \ref{modeling}: $D>0$ is the elastic stiffness parameter, while the $b_1 \in \mathbb R$ corresponds to in-axis forces and $b_2 \ge 0$ to the nonlinear restoring force; the function $p_0(x)$ represents static pressure on the surface of the beam, with flow parameters $\beta>0$ (density) and $U\ge 0$ (unperturbed flow velocity).  The parameter $k_0 \ge 0$ gives the {\em frictional} (or weak) damping coefficient, while  $k_1 \ge 0$ measures {\em square root} type damping (discussed in Section \ref{sqrt}). 

We distinguish between {\em panel}  and {\em cantilever} configurations below, with a key distinction being the allowance of rotational inertia. 
For the panel configurations, {\bf C} and {\bf H}, the dynamics are given by: 
\begin{equation}\label{Bergerplate}
\begin{cases} w_{tt}+D\partial_x^4 w + k_0w_t+[b_1-b_2||w_x||^2]w_{xx} =p_0(x)-  \beta [w_t+Uw_x] \\
BC(w) = \begin{cases} w(0)=w_x(0)=w(L)=w_x(L)=0 &~\text{for configuration  {\bf C } } \\[.1cm]
w(0)=w_{xx}(0)=w(L)=w_{xx}(L)=0 &~\text{for configuration  {\bf H } } \end{cases}\end{cases}.
\end{equation}
 For the cantilever configuration $\mathbf{CF}$, we allow $\alpha \ge 0$:
\begin{equation}\label{Bergerplate0}
\begin{cases}(1-\alpha\partial_x^2)w_{tt}+D\partial_x^4 w + [k_0-k_1\partial_x^2]w_t+[b_1-b_2||w_x||^2]w_{xx} =p_0(x)-  \beta [w_t+Uw_x] \\
 w(0)=w_x(0)=0 \\ 
  w_{xx}(L)=0, ~~ -\alpha \partial_xw_{tt}(L)+Dw_{xxx}(L)-k_1\partial_xw_t(L)+[b_1-b_2||w_x||^2]w_x(L)=0
 \end{cases} 
\end{equation}
In both cases above, appropriate initial conditions are specified: $w(t=0)=w_0;~~w_t(t=0)=w_1$.

{\em Our numerical studies below will compare dynamics across all three configurations holding {\bf $\alpha=0$}}; however, in the discussion of theoretical results, we will address technical issues associated with inertia and damping specific to the {\bf CF} configuration, and thus in said case we allow for $\alpha \ge 0$.

\begin{remark}
Though the clamped conditions have nice mathematical properties (e.g., $H_0^2$ displacements and preservation of regularity via extension by zero), the hinged boundary conditions are the easiest to work with numerically because the in vacuo modes are purely sinusoidal. Including a clamped or free condition gives rise to hyperbolic modes, as well as eigenvalues obtained via nontrivial transcendental equations. Additionally, the free end yields the aforementioned complications. On the other hand, we will see that it is easiest to destabilize the {\bf CF} configuration, followed by {\bf H}, with the {\bf C} beam most resistant to the flutter instability.
\end{remark}

\subsection{Rotational Inertia and Damping}\label{sqrt}
Standard aeroelasticity literature \cite{dowellnon,dowellchaos,vedeneev} omits the structural rotational inertia. A scaling argument is typically invoked for ``thin" structures, yielding $\alpha=0$ \cite{slenderness, lagnese}. Mathematically, the presence of $\alpha>0$ is regularizing and nontrivial: $w_t \in H^1$ \cite{springer}. Thus, $\alpha >0$ is helpful, and results with $\alpha=0$ are considered stronger \cite{webster,delay,springer,Memoires}. As such, the focus of the numerical study here is on $\alpha=0$. For the  {\bf CF} configuration the theoretical situation surrounding $\alpha$ is more complicated. The recent work \cite{HTW} fully explores well-posedness and long-time behavior (global attractors \cite{springer}) in the {\bf CF} configurations. We present these results below in Section \ref{sec:theory}.

To qualitatively study post-flutter dynamics, we here consider a simple nonlinear restoring force \cite{wonkrieg,beam4} described above. In the panel configurations {\bf C} and {\bf H}, the scalar nonlinearity can be treated as a locally Lipschitz perturbation, and the nonlinearity introduces no boundary terms. In the {\bf CF} configuration the nonlinearity {\em interacts} with the free boundary condition \eqref{badBC}, and non-trivial boundary traces appear in the analysis of trajectory differences. These boundary terms are only controlled with $\alpha>0$ (or by introducing boundary damping \cite{beam4} or other velocity smoothing).

Another interesting issue pertains to the strength of appropriate structural damping mechanisms {\em as a function of} $\alpha \ge 0$. The term $k_0$ measures the strength of standard {\em weak} structural damping. For $\alpha=0$, this is an appropriate strength of damping for stabilization; when $\alpha>0$, it is clear from the standard energy expression (see Section \ref{energysec}) that one requires $k_1>0$ in order for the damping to control of kinetic energies.

(Interesting spectral questions emerge for the less natural cases of $\alpha=0$ with $k_1>0$ and $k_0=0$, as well as $\alpha>0$ with $k_1=0$ and $k_0$ large---addressed below.)
When $k_1>0$, we say that the damping is of the {\em ``square root"} type, as $\partial_x^2$ is formally ``half'' the order of the principal stress operator $B=\partial_x^4$ \cite{beamdamping}. This concept can be generalized to fractional powers $[B]^\theta w_t$ for $\theta\in[0,1]$---Kelvin-Voigt damping occurs at $\theta=1$; see the more detailed discussion in \cite{HTW}, and references \cite{russell:93:JMAA,che-tri:89:PJM,han-las:00}.

Lastly, we point out the  clear disparity between the inclusion of inertia $\alpha>0$, and the aerodynamic damping provided from the flow itself in \eqref{linpist}: piston theory provides damping of the form $\beta w_t$ (with dissipative sign), yet, if one includes  inertia, this ``natural" weak damping is not adequate to stabilize energies. Thus, we have a mathematical incentive to take $\alpha=0$ if we wish to utilize the flow-provided damping to help bound or stabilize trajectories. 

\subsection{Notation}
For a given domain $D$,
its associated $L^{2}(D)$ will be denoted as $||\cdot ||_D$ (or simply $%
||\cdot||$ when the context is clear). Inner products in $L^{2}(D)$ are written $(\cdot ,\cdot)_{D}$ (or simply $(\cdot ,\cdot)$ when the context is clear). We will also denote pertinent duality pairings as $%
\left\langle \cdot ,\cdot \right\rangle _{X\times X^{\prime }}$, for a given
Hilbert space $X$. The space $H^{s}(D)$ will denote the Sobolev space of
order $s$, defined on a domain $D$, and $H_{0}^{s}(D)$ denotes the closure
of $C_{0}^{\infty }(D)$ in the $H^{s}(D)$-norm $\Vert \cdot \Vert
_{H^{s}(D)} $ or $\Vert \cdot \Vert _{s}$.

\subsection{Energies and State Space}\label{energysec}

 The natural energy for linear beam dynamics is given by the sum of the potential and kinetic energies
\begin{equation*}
E(t)=E_{\alpha}(w,w_t)=\frac{1}{2}\Big\{D||w_{xx}(t)||^{2}+
||w_{t}(t)||^{2}+\alpha ||\partial_xw_t||^2\Big\}.
\end{equation*}
Enforcing that $\alpha =0$ for {\bf H} and {\bf C} configurations, and that $\alpha \ge 0$ for {\bf CF}, we have that 
the dynamics evolve in the state space $\mathcal H$, whose definition depends on the configuration: 

\begin{equation}\label{statespace}\mathcal{H}  = \begin{cases} & \mathcal H_{\bf C} = H_0^2(0,L) \times L^2(0,L) \\[.1cm] &
\mathcal H_{\bf H}= (H^2\cap H_0^1)(0,L) \times L^2(0,L) \\[.2cm] &
\mathcal H_{\bf CF}= \begin{cases}  H^2_*(0,L)\times L^2(0,L) & \text{for}~\alpha=0 \\[.1cm] H^2_{*}(0,L)\times H^1(0,L) & \text{for}~\alpha>0\end{cases} \end{cases},\end{equation} where $$H_{*}^{2}(0,L)\equiv \{ w \in
H^{2}(0,L): w=w_x = 0 ~\text{ for }~x=0\}.$$ 

Owing to the conditional structure of the state space (depending on the configuration and, for {\bf CF}, the value of $\alpha$) we opt for a compact notation, using 
\begin{equation}\label{statespaceh}  H = \begin{cases}    H_0^2(0,L) & \text{ in configuration}~{\bf C} \\ 
 (H^2\cap H_0^1)(0,L) & \text{ in configuration}~{\bf H}   \\ 
H^2_{*}(0,L) & \text{ in configuration}~{\bf CF}  \end{cases},\end{equation}  with corresponding inner-product $(u,w)_H\equiv D(u_{xx},w_{xx})$ (equivalent to the standard $H^2(0,L)$ inner product), and $$L^2_{\alpha}\equiv \begin{cases} L^2(0,L) & \text{ for }~ \alpha=0 \\ H^1(0,L) & \text{ for}~\alpha>0\end{cases},~~~||w||^2_{L^2_{\alpha}}=\alpha||w_x||^2+||w||^2.$$ 
We can then cleanly write the state space across all cases as $$\mathcal H = H \times L^2_{\alpha},$$ with a clear meaning in a given context.
We
will also critically utilize the following nonlinear energies
associated to equations in (\ref{Bergerplate}) and \eqref{Bergerplate0}:
\begin{equation*}
\mathscr{E}(t)=\mathscr{E}(w,w_t)=E(t)+\Pi (w(t)),\hskip1cm
\widehat{E}(t)= E(t)+\frac{b_2}{4}||\nabla w(t) ||^{4},
\end{equation*}
where the $\Pi$ term represents the non-dissipative and nonlinear
portion of the energy:
\begin{equation}
\Pi(w)=\dfrac{1}{4}\Big(b_2||\nabla w||^{4}-2b_1||\nabla w||^2 \Big)-(p_0,w)
. \label{pi}
\end{equation}
\noindent As in the general case \cite[Lemma 1.5.4]{springer}, we can see that for any $w \in
 H$, $0<\eta\leq 2$ and
$\epsilon
>0,$
\begin{equation}\label{lot}
||w||_{2-\eta}^2 \le \epsilon \left[ ||w_{xx}(t)||^{2}+||w_x(t)||^{4}\right] +M(\epsilon).
\end{equation}
This yields the crucial fact that the positive nonlinear energy is dominant:
\begin{equation}
c_0\widehat{E}(w,w_t)-C\leq \mathscr{E}(w,w_t)\leq
c_1\widehat{E}(w,w_t)+C, \label{nonen}
\end{equation}%
for some $c_0,c_1,C>0$ depending on $b_1,b_2$. 

\section{Theory and Discussion}\label{sec:theory}
\subsection{Solutions}
We now discuss the pertinent notions of solution. We will distinguish between: (i) panel configurations ($\mathbf C$ and $\mathbf H$), for which we consider only $\alpha=0$, and (ii) the cantilever {\bf CF}, where we admit $\alpha \ge 0$. 
 \begin{itemize} \item {\em Weak solutions} satisfy a variational formulation of \eqref{Bergerplate}. One of the key features of such solutions is that $w_{tt}$ is interpreted only distributionally. 

   \item {\em Strong solutions} are weak solutions possessing additional regularity that permits a classical (point-wise) interpretation of the evolution \eqref{Bergerplate}. For {\bf CF} with $\alpha>0$, boundary conditions are still interpreted weakly due to subtleties associated with inertial terms. (See \cite{HTW} or \cite{springer}.)

   \item {\em Generalized solutions} are  $C([0,T],\mathcal H)$ limits of strong solutions (when they exist). These are weak solutions \cite{springer,Memoires}, but they admit smooth approximation, and thus inherit some properties of strong solutions holding with respect to the topology of $\mathcal H$. 
   \end{itemize}

\noindent Precise definitions are now given. For weak solutions, the variational form of the problem is the same for all three {\bf C}, {\bf H}, {\bf CF}, with the space $\mathcal H$ encapsulating the configuration and $\alpha \ge 0$.
   \begin{definition}[Weak Solutions]\label{def:weak0}
  We say  $w \in H^1(0,T; L^2(0,L))$ is a {\em weak} solution on $[0,T]$ if:
  \[ (w,w_t) \in L^{\infty}(0,T;\mathcal H)\cap C_w([0,T]; \mathcal H)~~\] 
and, for every $\phi \in H$, we have
  \begin{align}\label{weak-sol0}
       \frac{d}{dt} \big[(w_t, \phi)+ \alpha (w_{tx}, \phi_x)  \big] + (w,\phi)_{H} + (b_2 \|w_x\|^2-b_1)(w_x,\phi_x) 
      + &k_0(w_t,\phi) + k_1 (w_{tx},\phi_x)\\ \nonumber
  & = (p_0, \phi) -\beta (w_t, \phi) -\beta U (w_x,\phi),
  \end{align}
  where $d/dt$ is taken in the sense of $\cD'(0,T)$. Moreover, for any $(\chi,\psi) \in \mathcal H$
  \begin{equation}\label{weak-ic0}
    (w,\chi)_{H}\big|_{t\to 0^+} = (w_0, \chi)_{H},\quad (w_t,\psi)\big|_{t\to 0^+} = (w_1,\psi).
  \end{equation}
\end{definition}
\begin{definition}[Panel strong solution]
For {\bf C} or {\bf H} configurations, a weak solution to \eqref{Bergerplate} is {\em strong} if it possesses the regularity: $(w,w_t)\in L^\infty(0,T; H^4(0,L)\times H)$ and $\frac{d^+}{dt^+}w_t$ is right-continuous. 
\end{definition}

\noindent More regular cantilever solutions are considered, but only for $\alpha>0$. 
\begin{definition}[Cantilever strong solution] For the {\bf CF} configuration,
  we say a weak solution to \eqref{Bergerplate0} with $\alpha>0$ is strong if it possesses the regularity: $w\in L^\infty(0,T; W)$, where 
\begin{equation}\label{def:W}
  W \equiv \{ v \in (H^2_*\cap H^3)(0,L) : v_{xx}(L) = 0\},
\end{equation}
    and $w_t\in L^\infty(0,T; H^2_*(0,L))$. In addition $\frac{d^+}{dt^+}w_t$ is right-continuous and $L^\infty(0,T;H^1_*(0,L))$ with $H^1_*(0,L) \equiv \{w \in H^1(0,L)~:~w(0)=0\}$.  (Such solutions satisfy the third-order boundary condition weakly---see \cite{HTW}.)
\end{definition}
 \noindent Finally, the most convenient notion of solutions will be those induced by a semigroup flow:
\begin{definition}[Semigroup well-posedness]\label{def:semigroup}
  We say \eqref{Bergerplate} or \eqref{Bergerplate0} is semigroup well-posed if there is a family of locally Lipschitz operators $t\mapsto \mathcal S(t)$ on $\mathcal H$, such that: for any $y_0=(w_0,w_1)\in \mathcal H$ the function $t\mapsto S(t)y_0$ (i) is in $C([0,T]; \mathcal H)$, (ii) is a weak solution to \eqref{Bergerplate}, and (iii) is a strong $C([0,T];\mathscr H_{\alpha})$ limit of strong solutions to \eqref{Bergerplate} on every $[0,T]$, $T >0$. In particular, solutions are unique and depend continuously in $C([0,T];\mathcal H)$ on the initial data from $\mathcal H$. Thus $t\mapsto S(t)y_0$ is a {\em generalized solution}.
   Furthermore, there exists some subset of $\mathcal H$, denoted $\mathscr D(\mathbb A)$ (with $\mathbb A$ the nonlinear generator of $\mathcal S(t)$), invariant under the flow, such that all solutions originating there are {\em strong solutions}.
  \end{definition}

\subsection{Well-posedness}
For all results below $D,k_0,\beta,U \ge 0$. We employ the energetic notations from Section \ref{energysec}.
\begin{theorem}[Linear semigroup well-posedness] For configurations {\bf C} and {\bf H}, the linear problem in \eqref{Bergerplate} (with $b_1=b_2=0$) is semigroup well-posed on $\mathcal H_{\mathbf C}$ or $\mathcal H_{\mathbf H}$ (resp.). 
Generalized (and strong) solutions obey the energy identity for all $t \ge 0$:
\begin{equation}
E(t) +(k_0+\beta)\int_0^t ||w_t(\tau)||^2d\tau = E(0) +\int_0^t(p_0,w_t)d\tau- \beta U \int_0^t \int_0^L [w_x(\tau)][w_t(\tau)] dx d\tau
\end{equation}
Similarly, for  configuration {\bf CF} with $\alpha, k_1 \ge 0$, the linear problem in \eqref{Bergerplate0} (with $b_1=b_2=0$) is semigroup well-posed on $\mathcal H_{\mathbf{CF}}$. 
Generalized (and strong) solutions obey the energy identity for all $t \ge 0$:
\begin{equation}
E_{\alpha}(t) +\int_0^t\Big[(k_0+\beta) ||w_t(\tau)||^2+k_1||\partial_xw_t||^2\Big]d\tau = E_{\alpha}(0)+\int_0^t(p_0,w_t)d\tau - \beta U \int_0^t \int_0^L [w_x(\tau)][w_t(\tau)] dx d\tau
\end{equation}
\end{theorem}
\noindent For such a linear problem, these results are well-established. 

\begin{theorem}[Nonlinear semigroup well-posedness]
Let $b_1 \in \mathbb R$ and $b_2 >0$.  For configurations {\bf C} and {\bf H}, the problem in \eqref{Bergerplate}  is semigroup well-posed on $\mathcal H_{\mathbf C}$ or $\mathcal H_{\mathbf H}$ (resp.). 
Generalized (and strong) solutions obey the energy identity for all $t \ge 0$:
\begin{equation}\label{energy1}
\mathscr E(t) +(k_0+\beta)\int_0^t ||w_t(\tau)||^2d\tau = \mathscr E(0) - \beta U \int_0^t \int_0^L [w_x(\tau)][w_t(\tau)] dx d\tau
\end{equation}
Similarly, for  configuration {\bf CF} with $\alpha> 0$ and $k_1 \ge 0$, the problem in \eqref{Bergerplate0} is (nonlinear) semigroup well-posed on $\mathcal H_{\mathbf{CF}}$. 
Generalized (and strong) solutions obey the energy identity for all $t \ge 0$:
\begin{equation}\label{energy2}
\mathscr E_{\alpha}(t) +\int_0^t\Big[(k_0+\beta) ||w_t(\tau)||^2+k_1||\partial_xw_t||^2\Big]d\tau = \mathscr E_{\alpha}(0) - \beta U \int_0^t \int_0^L [w_x(\tau)][w_t(\tau)] dx d\tau
\end{equation}
\end{theorem}
In the case of nonlinear clamped and hinged beams, the well-posedness has been known for sometime; consult the abstract theory in \cite{springer} and note the classical references (for both 1-D Krieger and 2-D Berger nonlinearities) \cite{ball, ball2,biacrip,dickey, edenmil}. 
In the case of the cantilevered beam {\bf CF} with $\alpha \ge 0$, see the recent work \cite{HTW} (which includes discussion of other cantilevered beams \cite{beam4}). 

The final case not addressed by the previous theorem is the nonlinear cantilever in the absence of rotational inertia. For this case \cite{HTW} provides only an existence result.
\begin{theorem}\label{weaksols} Let $b_1 \in \mathbb R$ with $b_2 \ge 0$ and $k_1 \ge 0$.
For \eqref{Bergerplate0} with $\alpha=0$ in the {\bf CF} configuration, {\em weak solutions} exist on $[0,T]$ for any $T>0$. 
\end{theorem}
Without boundary damping (e.g., \cite{beam4}), existence of strong solutions for {\bf CF} with $\alpha=0$ is open.

\subsection{Long-time Behavior of Trajectories}\label{long-time}
We can obtain the following global-in-time stability bounds for all generalized solution to \eqref{Bergerplate} or \eqref{Bergerplate0} in the presence of piston-theoretic terms, $\beta \ge 0$ and $U \ge 0$:
\begin{theorem} Let $D >0$ with $\beta>0$ and $U \ge 0$.
In configurations {\bf C} and {\bf H}, we have the following global in time bound for generalized solutions to \eqref{Bergerplate} for any amount of weak damping $k_0 \ge 0$:
\begin{equation}\label{energybound1}\sup_{t \ge 0} \Big\{||w_{xx}(t)||^2+||w_t||^2\Big\} \le C.\end{equation}
In the case of configuration {\bf CF}, {\em with $\alpha>0$} (rotational inertia present) and $k_1>0$ (active square root damping), we have 
\begin{equation}\label{energybound2}\sup_{t \ge 0} \Big\{||w_{xx}(t)||^2+||w_t||^2+\alpha||\partial_xw_t||^2\Big\} \le C.\end{equation}
\end{theorem}

\begin{remark}
Formally, {\bf CF} trajectories for \eqref{Bergerplate0} with no inertia, $\alpha=0$, do obey a global-in-time bound of the form \eqref{energybound1}, but the requisite multipliers (and associated integration by parts) are not permissible for the weak solutions guaranteed by Theorem \ref{weaksols}.
\end{remark}

These global bounds are nontrivial, owing to the fact that the dynamics are non-gradient, as reflected by the energy-building terms under temporal integration in \eqref{energy1} and \eqref{energy2}. 
They follow from the estimates in \eqref{lot}--\eqref{nonen}, taken together with a Lyapunov argument, and rely critically on the ``good" structure of $\Pi$, reflected by the control of lower order terms in \eqref{lot}. In fact, more is true. For a dynamical system, a {\em compact global attractor} is a uniformly attracting (in the sense of the $\mathcal H$ topology), fully time-invariant set \cite{springer}. 
\begin{theorem} Let $D>0$ with $\beta>0$ and $U\ge 0$. Also let $b_2>0$ with $b_1 \in \mathbb R$.

For configurations {\bf C} and {\bf H} with $k_0 \ge 0$ the dynamical system $(S(t),\mathcal H)$ generated by generalized solution to \eqref{Bergerplate} has a smooth and finite dimensional compact global attractor. 

For configurations {\bf CF} with $\alpha,k_1 > 0$ and $k_0 \ge 0$, the dynamical system $(S_{\alpha}(t),\mathcal H_{\mathbf{CF}})$ generated by generalized solution to \eqref{Bergerplate0} has a finite dimensional compact global attractor. 
\end{theorem}
In the above theorem, {\em smooth} is taken to mean that trajectories on the attractor have $(w,w_t)\in \mathcal H \cap \big(H^4(0,L) \times H\big)$, where $H$ indicates the configuration-dependent displacement space from \eqref{statespaceh}; the attractor is bounded in this topology. {\em Finite dimensional} refers the fractal dimension of the attractor---see \cite{springer} for more details.

In this study, the attractor ``contains" the essential post-flutter dynamics for 
\eqref{Bergerplate} or \eqref{Bergerplate0}. The aeroelastic model is non-conservative, thus the associated dynamical system has no strict Lyapunov function \cite{springer}. In this case the global attractor is not characterized as the unstable manifold of the equilibria set, so to obtain a compact global attractor one must show both {\em asymptotic compactness} and {\em ultimately dissipativity} of the dynamical system \cite{springer}.  Since the infinite dimensional dynamical system associated to solutions has the property that trajectories converge---in a uniform sense---to an \textit{attractor}, we effectively reduce the analysis of the dynamics to a ``nice" finite dimensional set.

The existence and properties of the attractor for panel configurations are shown in \cite{HLW} (though earlier proofs have certainly appeared for these configurations). This work also discusses the application of a recent tool: the theory of {\em quasi-stability} \cite{quasi,springer}. The theory provides the existence of so called {\em exponential attractors} \cite{quasi,springer}. 
A {\em generalized fractal exponential attractor} for the dynamics $(S(t), \mathcal H)$ is a forward invariant compact set $A_{\text{exp}} \subset \mathcal H$, with finite fractal dimension (possibly in a weaker topology), that attracts bounded sets with uniform exponential rate.  We remove the word ``generalized" if the fractal dimension of $A_{\text{exp}}$ is finite in $\mathcal H$. 

\begin{theorem} Let $D>0$ with $\beta>0$ and $U\ge 0$. Also let $b_2>0$ with $b_1 \in \mathbb R$.

For configurations {\bf C} and {\bf H} with $k_0 \ge 0$ the dynamical system $(S(t),\mathcal H)$ generated by generalized solution to \eqref{Bergerplate} has a generalized fractal exponential attractor. 

If the imposed damping is sufficiently large, i.e.,~$k_0>k_*(U,\beta,b_i,D,L,\beta)>0$, the dynamical system $(S(t),\mathcal H)$ generated by generalized solutions to \eqref{Bergerplate} has a fractal exponential attractor. 
\end{theorem}
\noindent In the cantilever configuration {\bf CF}, the result on the existence of a compact global attractor is much more recent  \cite{HTW}.

\section{Overview of Numerical Studies}\label{mainresults}

Our numerical investigations of the piston-theoretic beam are multi-faceted; we seek to determine critical parameter values that lead to the onset of instability, and  to understand essential qualitative properties of post-flutter dynamics. To be consistent with comparable engineering studies, and to be consistent across configurations, we neglect rotational inertia here.
\begin{center} {\bf In all of our experiments, we focus on the case $\alpha=0$ (and thus $k_1=0$).}\end{center} We also take the static pressure $p_0(x) \equiv 0$.  In Section \ref{linear} we use physical parameters to emphasize the predictive capabilities of our methods. In Section \ref{postflutter}, our primary interest is the qualitative properties of post-flutter nonlinear dynamics, and, as such, we choose mathematically convenient parameter values.

We utilize two main methods for these investigations: {\em modal analysis} (described in detail in Section \ref{modal}) and {\em finite differences}. To analyze the onset of instabilities in the linear model, we use modal analysis as a predictive tool for finding critical parameters, and we use finite difference approximations to determine if these predictions are in agreement with the dynamics. Subsequently in Section \ref{postflutter}, we peform finite difference simulations of the nonlinear model to investigate the post-flutter regime.  

Finite difference simulations are accomplished via standard second-order difference formulas for spatial derivatives, with ghost values beyond the spatial boundary as necessary for the various boundary configurations.  Temporal integration is accomplished using a stiff ODE solver.  Our simulations utilize three basic types of initial data (configuration dependent): (i) {\em modal initial displacement}, (ii) {\em polynomial initial displacement}, and {\em elementary initial velocity}. Let us denote $\hat x=x/L$.
\begin{itemize}
\item ~\noindent \textbf{[$n$th Mode ID]}  For (i), we  consider the  in-vacuo mode shapes of a given  configuration (detailed below in Section \ref{modal}) as an initial displacement, with zero initial velocity profile. We restrict to mode numbers $n=1,..., 5$. For instance, in the case of {\bf CF} conditions we have:
$$w(0,\hat x)=s_2(\hat x) =  [\cos(k_2\hat x)-\cosh(k_2\hat x)]-\mathcal C_2[\sin(k_2\hat x)-\sinh(k_2\hat x)], ~~~~~w_t(0,\hat x)=0,$$  
where $k_2\approx 4.6941$ is the second Euler-Bernoulli cantilevered mode number and $\mathcal C_2 \approx 1.0185$.
 
\item ~\textbf{[Polynomial ID]} For (ii), with the {\bf C} and {\bf F} configurations, we  consider a polynomial of minimal degree {\em satisfying the boundary conditions} for the respective configuration in the experiment.  For {\bf CF,} we have $$w(0,\hat x)=-\hat 4x^5+15\hat x^4-20\hat x^3+10\hat x^2, ~~~~~w_t(0,\hat x)=0.$$
We  again take the initial velocity profile to be zero.
 
\item ~\textbf{[Elementary IV]} For (iii), we take the initial displacement to be identically zero, and consider a simple initial velocity profile. In the case of {\bf C} or {\bf H}, we will take a symmetric quadratic function on $[0,L]$; for {\bf CF} we simply take a linear function of positive slope on $[0,L]$ vanishing at $x=0$: 
  $$w(0,\hat x)=0,~~~~~w_t(0,\hat x)=\hat x.$$
\end{itemize}

In the sections to follow, we present our numerical results.  In some sections we compare across boundary configurations, but for qualitative studies we often choose a single, particularly illustrative example in a given configuration.  In addition, we sometimes consider how qualitative aspects of the simulations may or may not vary across differing initial conditions.  Several visualizations involve energy plots (as a function of time), as well as plots of midpoint displacements (in time) for panel ({\bf C} or {\bf H}) configurations, or free end displacements for the cantilever ({\bf CF}) configuration.

In Section \ref{end} we provide a synopsis of the main numerical observations below, and connect to the theory presented in the previous sections.

\section{Linear Computational Analysis: Flutter Prediction}\label{linear}
In this section we analyze the {\em onset} problem for flow-induced instability: namely, what combinations of parameters of interest---$U,\beta$ for the flow, $k_0,L$ for the beam---yield stable or unstable asymptotic-in-time dynamics. More specifically, our studies in this section we are concerned with {\em convergence to stationary states} or {\em unbounded growth of displacements} as $t\to \infty$. (Note: when $b_2=0$, without the nonlinear restoring force, unstable trajectories will grow unboundedly.) 

For the readers convenience, we restate the system studied in this section here:
\begin{equation}\label{Bergerplate**}
\begin{cases} w_{tt}+D\partial_x^4 w + k_0w_t = -  \beta [w_t+Uw_x] \\[.1cm]
w(t=0)=w_0;~~w_t(t=0)=w_1\\[.2cm]
BC(w) = \begin{cases} w(0)=w_x(0)=w(L)=w_x(L)=0 &~\text{for configuration  {\bf C } } \\[.1cm]
w(0)=w_{xx}(0)=w(L)=w_{xx}(L)=0 &~\text{for configuration  {\bf H } } \\ 
w(0)=w_x(0)=w_{xx}(L)=w_{xxx}(L)=0& ~\text{for configuration} ~{\bf CF} \end{cases}\end{cases}.
\end{equation}
For convenience, let us define the damping parameter ~$k\equiv k_0+\beta$,~ which consists of both imposed damping $k_0 \ge 0$, as well as piston-theoretic (flow) damping $\beta >0$.

The effects of the initial conditions are conservative (which is to say that, for $\beta=k_0=0$, the energy is conserved throughout deflection).
With no imposed damping ($k_0=0$), the transient dynamics may or may not damp out; in this case, the only damping in the problem, then, is contributed by the flow in the form of $\beta w_t$. Additionally, with $U>0$, the non-dissipative piston term $-\beta Uw_x$ on the RHS of \eqref{Bergerplate} may induce instability (and thus oscillatory behavior in the linear energy $E(t)$). In this sense, stability can be seen as competition between the effect of the damping in $\beta w_t$ bleeding energy out of the system, and the non-dissipative flow effect $Uw_x$ putting flow energy in. As $U$ is increased the effect of the flow becomes more pronounced and the location of the spectrum is shifted. For all other parameters fixed (including $\beta$), a certain value $U_{\text{crit}}$ represents a bifurcation, and the dynamics become unstable and exponential growth (in time) of the solution is observed. 
With imposed damping in the model ($k_0>0$ and $L$, $\beta$ fixed), there exists a $U_{\text{crit}}(k_0)$ such that for $U<U_{\text{crit}}$, the transient dynamics are damped out, and for $U>U_{\text{crit}}$ instability will still be observed. We can think of $U_{\text{crit}}$ as an increasing function of the damping $k_0$ in the problem.  

\subsection{Modal Analysis}\label{modal}
Modal analysis, often referred to in engineering literature, can in fact have various meanings. Generally, it refers to a Galerkin method whereby solutions to a structural or fluid-structural problem are approximated by in vacuo structural eigenfunctions (``modes") (e.g., \cite{modalcant,jfs,vedeneev,dowellnon,inext1}). Since the eigenfunctions of standard elasticity operators for a basis for the state space,  a good well-posedness result for the full system justifies this type of approximation. This type of approximation can be dynamic, as in reducing an evolutionary PDE to a finite dimensional system of ODEs by truncation, or it can be stationary, reducing the problem of dynamic instability (for linear dynamics) to an algebraic equation. The latter is what we utilize and describe in detail.

\subsubsection{Spatial Modes and Eigenvalues}
Critical to any modal analysis are the in vacuo modes (eigenfunctions) associated with each configuration: {\bf C}, {\bf H}, {\bf CF}. Since we are working with the Euler-Bernoulli  beam ($\alpha=0$) for our simulations, the modes and associated eigenvalues can be computed in an elementary way. These functions are {\em complete} and {\em orthonormal} in $L^2(0,L)$, as well as complete and orthogonal in $H^2(0,L)$ (with respect to $(\cdot,\cdot)_H$, where $H$ encodes the boundary conditions \eqref{statespaceh}).

For the in vacuo dynamics, we have:~~$w_{tt}+D\partial_x^4 w  = 0.$
Separating variables yields the dispersion relation: $\kappa^4_n = \dfrac{\omega_n^2}{D},$ where we have labeled temporal frequencies $\omega_n$ and associated mode numbers $\kappa_n$, $n=1, 2,....$ The spatial problem is then
$\partial_x^4s+\kappa_n^4 s=0$, giving the fundamental set\footnote{This choice of linear combinations is particularly convenient.}: 
 \begin{equation*}
\Big\{\cos(\kappa_nx)+\cosh(\kappa_nx),~~\cos(\kappa_nx)-\cosh(\kappa_nx),~~\sin(\kappa_nx)+\sinh(\kappa_nx),~~\sin(\kappa_nx)-\sinh(\kappa_nx)\Big\}.
 \end{equation*}
 So the modes  in any configuration are linear combinations of the above, determined by invoking the boundary conditions; the boundary conditions also
 provide the specific values $(\kappa_n,\omega_n)$.  
\begin{table}[h]
\begin{center}
\begin{tabular}{llll}
Config. & Mode Shape & Char. Equation & $\kappa_n$\\\hline 
{\bf H} & $\sin(\kappa_nx)+C_n\sinh(\kappa_nx)$ & $\sin(\kappa_nL)\sinh(\kappa_nL)=0$ & $\pi n$ \\
{\bf C} & $(\cos(\kappa_nx)-\cosh(\kappa_nx))+ C_n(\sin(\kappa_nx)-\sinh(\kappa_nx))$ & $\cos(\kappa_nL)\cosh(\kappa_nL)=1$ & Table \ref{table1}\\
{\bf CF}  & $(\cos(\kappa_nx)-\cosh(\kappa_nx))+ C_n(\sin(\kappa_nx)-\sinh(\kappa_nx))$  & $\cos(\kappa_nL)\cosh(\kappa_nL)=-1$ & Table \ref{table1}\\
\end{tabular}
\caption{Mode shapes and characteristic values for the various boundary configurations. In each case $C_n$ is obtained from the characteristic equation. }
\label{tab:t1}
\end{center}
\end{table}

For {\bf C} and {\bf CF} in Table \ref{tab:t1}, the mode numbers $\kappa_nL$ are obtained by numerically solving the characteristic equation.  We have
\[
{\bf C}: C_n=\dfrac{-c_n\big(\cos(\kappa_nL)-\cosh(\kappa_nL)\big)}{\sin(\kappa_nL)-\sinh(\kappa_nL)} \qquad\qquad {\bf CF}: C_n=\dfrac{-c_n\big(\cos(\kappa_nL)+\cosh(\kappa_nL)\big)}{\sin(\kappa_nL)+\sinh(\kappa_nL)}.\]
The $c_n$ values are chosen to normalize the  functions in the $L^2(0,L)$ sense.

\subsubsection{Reduction to Perturbed Eigenvalue Problem}\label{modalapproach}
Let us consider the Galerkin procedure for the full linear equation beam equation: \begin{equation}\label{fornumerics}w_{tt}+Dw_{xxxx}+kw_t=-\beta Uw_x\end{equation} on $(0,L)$, with boundary conditions according to the configuration being studied; recall that $k=k_0+\beta$. We view the terms involving $\beta, U,k_0$ as {\em perturbations} and expand the solution  via the in vacuo  mode functions $\{s_j\}$ as
$w(t,x) = \sum q_j(t)s_j(x)$.  The $\{q_j\}$ represent smooth, {\em time-dependent coefficients}. Plugging this representation into the \eqref{fornumerics}, multiplying by $s_n$, and integrating over $(0,L)$ for each $n$ we obtain:
\begin{align}
\sum_m\Big[\big[q_{m,tt}(t)+(\beta+k_0)q_{m,t}(t) +Dk_m^4q_m(t)\big](s_m,s_n)_{(0,L)}+\beta U(\partial_xs_m,s_n)_{(0,L)}q_m(t)\Big]=&~0,
\end{align}
Orthonormality of the eigenfunctions can be invoked to produce {\em diagonal} terms, whereas the terms scaled by $\beta U(\partial_xs_m,s_n)$ are {\em off-diagonal} and give rise to the instability of the ODE system. 
\begin{remark} To obtain a spectral simulation for the dynamics, one could simply truncate the system above for all $1<m,n<N$ and solve the resulting $N\times N$ system for the unknown coefficients $\{q_j(t)\}$, yielding the approximate solution $\displaystyle w_N(x,t)=\sum_{j=1}^Nq_j(t)s_j(x).$ \end{remark}

To simply determine the stability of the problem as a function of the given parameters, we invoke a standard engineering ansatz \cite{vedeneev,dowell} (and references therein): {\em assume simple harmonic motion in time of the beam} according to some dominant (perturbed) frequency $\widetilde\omega$; we allow possible contribution from all functions $s_n$ for $n=1,2,...,N$ via coefficients labeled $\alpha_n$:
\begin{equation}\label{preplug}w(t,x) \approx e^{-i\widetilde{\omega}t}\sum_{j=1}^{N}\alpha_j s_j(x),\end{equation} 
where $N$ is a chosen dimensional {\em truncation}. 
We repeat by multiplying by mode functions $s_n$, and then integration produces an eigenvalue problem in the perturbed frequency $\widetilde{\omega}$. 
With the off-diagonal entries $ (\partial_x s_m,s_n)$  in hand (for $1 \le m,n \le N$ with $m \neq n$), compute diagonal terms $$\Omega_j(\widetilde \omega) = -\widetilde \omega^2-i(\beta+k_0)\widetilde\omega +Dk_j^4,~~j=1,2,...,N,$$ we create the matrix for $1 \le n,m \le N$: 
\begin{equation}
A =A(\widetilde \omega)=[a_{mn}],~~\text{ with}~~a_{mn} = \begin{cases} \Omega_m &~\text{ for }~m =n \\
\beta U (\partial_xs_n,s_m) &~\text{ for }~m \neq n.
\end{cases}\label{modalmatrix}
\end{equation}

As an example, for the configuration {\bf H}, the entries in the matrix can be computed explicitly. Setting $N=4$ yields the linear system
\begin{equation} \begin{bmatrix}
\Omega_1 &  \dfrac{-8 \beta U }{3L} & 0  &\dfrac{-16 \beta U }{15L} \\[.3cm]
\dfrac{8 \beta U }{3L} & \Omega_2 &\dfrac{-24 \beta U }{5L}&0 \\[.3cm]
0 & \dfrac{24 \beta U }{5L}& \Omega_3 &\dfrac{-48 \beta U }{7L}  \\[.3cm]
 \dfrac{16 \beta U }{15L} &0&\dfrac{48 \beta U }{7L} & \Omega_4
 \end{bmatrix} \begin{bmatrix} \alpha_1 \\[.5cm] \alpha_2  \\[.5cm] \alpha_3  \\[.5cm] \alpha_4   \end{bmatrix} = \mathbf 0.
\end{equation}
For chosen parameter values of $D, k_0, \beta, U, L$, we enforce the zero determinant condition for non-trivial solutions in $\widetilde{\omega}$:
$$\text{det}\big(A(\widetilde \omega)\big) =0,$$ and solve for $\widetilde \omega = [\widetilde \omega_{1},...,\widetilde \omega_N]^T$. The associated complex roots allows us to track the response of the natural modes to the {\em perturbation} terms---damping, $(k_0+\beta)w_t$ and non-dissipative, $\beta Uw_x$.  
 \begin{itemize} \item If $Im(\widetilde{\omega}_j)>0$, then we call the configuration {\em unstable}, and the associated linear dynamics should exhibit exponential growth in time with rate $Im(\widetilde{\omega}_j)$ and frequency $Re(\widetilde{\omega}_j)$. 
\item Otherwise, the linear dynamics  remain bounded in time. \end{itemize}

Obviously, this method of determining stability has associated error---first, owing to the simple harmonic motion ansatz, and secondly, because of truncation.
Testing the conclusions drawn from this type of modal stability determination against a full numerical scheme is one of the main numerical points of this consideration.
\begin{remark}
The emergence of an unstable perturbed eigenvalue (associated to $Im(\widetilde \omega)$) is precisely what is meant by {\em self-excitation instability}, which is commonly called {\em flutter}. In this way, we note that predicting {\em the onset of flutter is linear}.  As mentioned above, for one to ``observe" (numerically) the post-onset/post-flutter dynamics, a nonlinear elastic restoring force must be incorporated into the model to keep solutions bounded in time (as the underlying linear dynamics are driven by the aforementioned exponential growth). \end{remark}

\begin{remark}
From a practical point of view, one can make use of any reasonable $N$ for modal calculations involving the {\bf H} configuration---these mode functions are purely sinusoidal. In modal calculations involving {\bf CF} or {\bf C}, we must first solve for the $k_n$, which is more numerically unstable for larger $n$.  Additionally, these approximate $k_n$ are used in the numerical integration of $\cosh$ and $\sinh$ functions with large arguments, further propagating error.
 Therefore, for $N\ge 7$ for {\bf C} and {\bf CF}, modal calculations are cumbersome and possibly inaccurate. Engineers have devised clever approximation of the mode functions is needed to circumvent this issue \cite{dowellapprox1}, though it does not concern us too much here as flutter is considered a ``low-mode" phenomenon. As evidenced below, modal approximations perform quite well, even for $N=6$. 
\end{remark}

\subsection{Prediction of the Onset of Instabilities}\label{predict}
\def\Ucrit{U_{\text{crit}}}
\def\hUcrit{\hat{U}_{\text{crit}}}
\def\emax{\mathcal{E}_{\text{max}}}

We now predict the onset of instability for the model \eqref{Bergerplate0} as the $U$ parameter increases, with other parameter values fixed and the nonlinear parameter $b_2=0$.  Indeed, there exists some critical flow speed value, $\Ucrit$, beyond which the system exhibits exponential growth in time.

In this section, we have adopted the non-dimensionalization procedure utilized in \cite{jfs}.  The parameters chosen here represent physical values for the prediction of the flutter phenomenon.  Accordingly, we utilize these as the fixed baseline values throughout this section:
 \begin{itemize}
 \item $D=23.9$ (fixed across all simulations),
  \item $L=300$ (allowed to vary in the range $100 \le L \le  500$,
 \item $\beta =1.2\cdot 10^{-4}$ (allowed to vary in the range $10^{-5} \le \beta \le 10^{-2}$).
 \end{itemize}

When employing finite difference simulations of \eqref{Bergerplate0}, the task of identifying $\Ucrit$ usually requires several batches of simulations, usually in some bracketing approach. In general it is not practical, as qualitative observations of flutter are somewhat subjective and may require excessive long-time simulations (for instance, accurately determining $\Ucrit$ for $L>320$).  However, the modal approach outlined in Section \ref{modalapproach} yields a direct, simple way to find $\Ucrit$. For given parameter values, one can immediately analyze the perturbed eigenvalues associated with \eqref{modalmatrix} to determine if an instability is present ($Im(\tilde{\omega}_j)>0$ for some $j$).  
 
\begin{figure}[H]
\begin{center}
\includegraphics[width=6in]{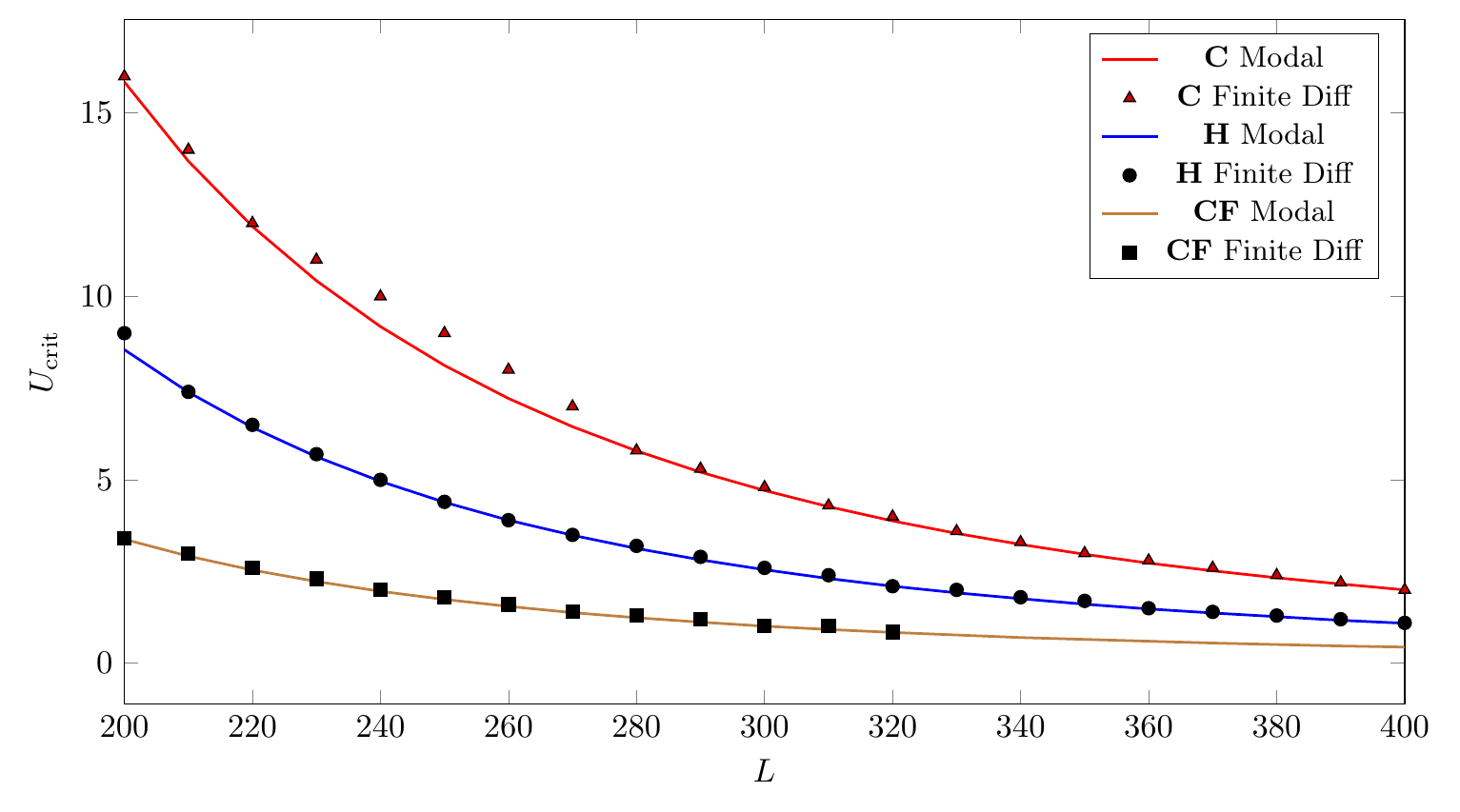}
\caption{Plot of critical $U$ values for $L$, $D=23.9$ and $\beta=1.2\cdot 10^{-4}$, $b_1=b_2=0$.}\label{fig1}
\end{center}
\end{figure}

In Figure \ref{fig1}, we compare the results of a bracketing search for $\Ucrit$ with the modal approach for a range of $L$ values across all three boundary configurations.  The points are the $\Ucrit$ values found by a bisection-like search employing finite difference simulations, and the curves are the values predicted by the modal approach. We see that the modal approach is very consistent with the finite difference approach, and is substantially more convenient, for predicting $\Ucrit$. The hierarchy of stability between the three configurations {\bf C}, {\bf H}, and {\bf CF} is also revealed in Figure \ref{fig1}

Note that the modal approach also produces an $\omega_{\text{crit}}$
associated to each $\Ucrit$ (not shown here), and the frequency of unstable oscillation is simply $Re(\omega_{\text{crit}})$.  Finally, as the modal analysis was such a good predictor of $\Ucrit$ for all configurations, the laborious finite difference simulations to search for $\Ucrit$ in the clamped-free case for $L>320$ were not performed.

\begin{figure}[H]
\begin{center}
\includegraphics[width=6in]{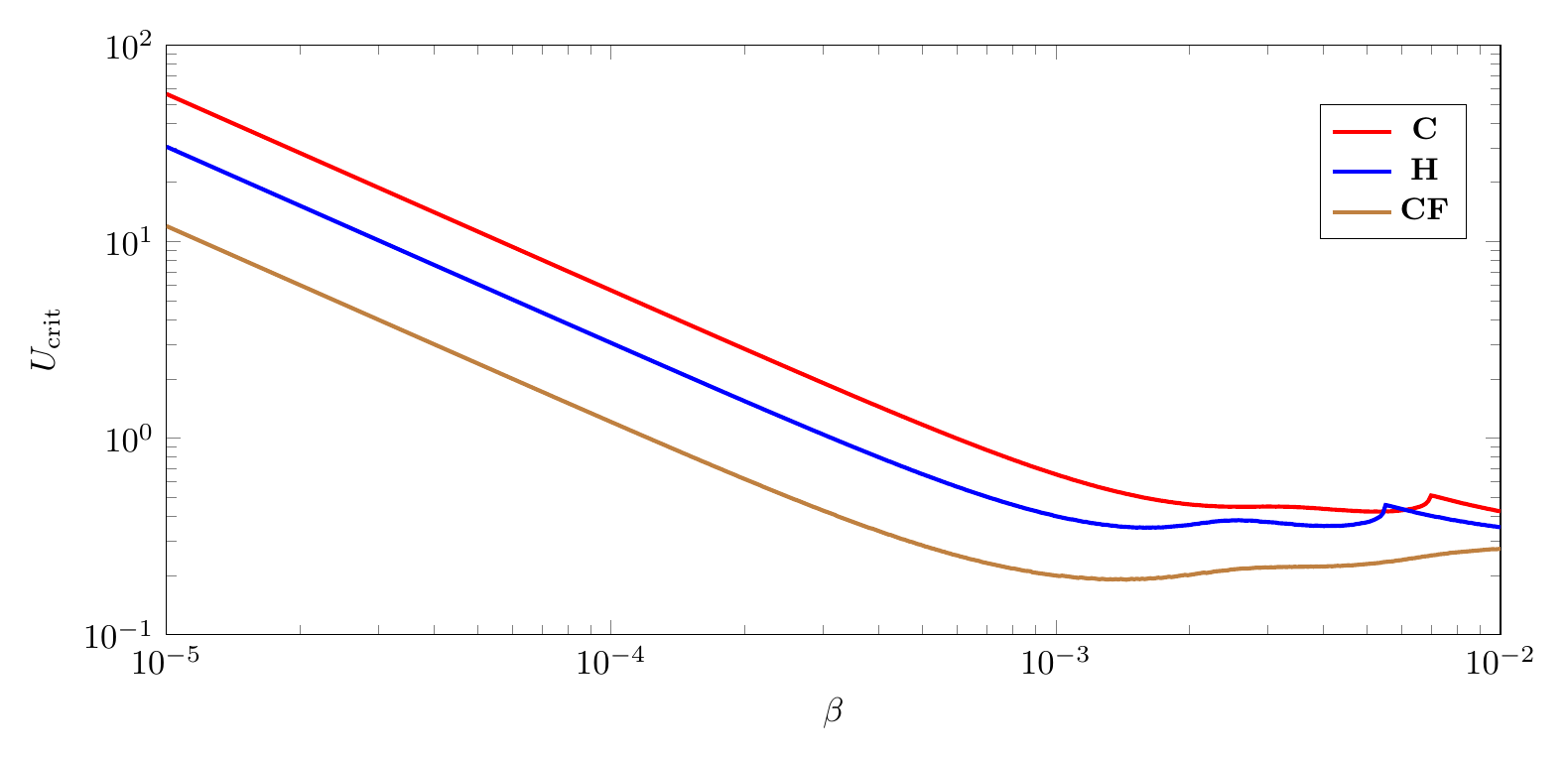}
\caption{Plot of critical $U$ values for $k_0$, $k=k_0+\beta$, $\beta=1.2\cdot10^{-4}$, $L=300$, $D=23.9$.}\label{fig2}
\end{center}
\end{figure}

In Figure \ref{fig2}, we determine $\Ucrit$ as a function of $\beta$, with  $10^{-5}\le \beta\le 10^{-2}$.  We observe monotonic linear decay in the  $\Ucrit(\beta)$ value as $\beta$ increases up to around $10^{-3}$, and then somewhat erratic behavior for larger $\beta$. Since $\beta$ affects both the damping and destabilizing terms in the model \eqref{fornumerics}, it is not obvious that its effects will remain monotonic on $\Ucrit$ as $\beta$ increases further.  Indeed, an analysis of the dispersion relation associated to \eqref{fornumerics} indicates that a critical transition occurs at $\beta=1$.

One could similarly extract critical values of $L$ or $\beta$ by fixing the other and varying $U$.  Results maintain the same qualitative structure, for instance, with respect to the established stability hierarchy of boundary configurations and monotonic behavior of $\Ucrit$. 

\begin{figure}[H]
\begin{center}
\includegraphics[width=6in]{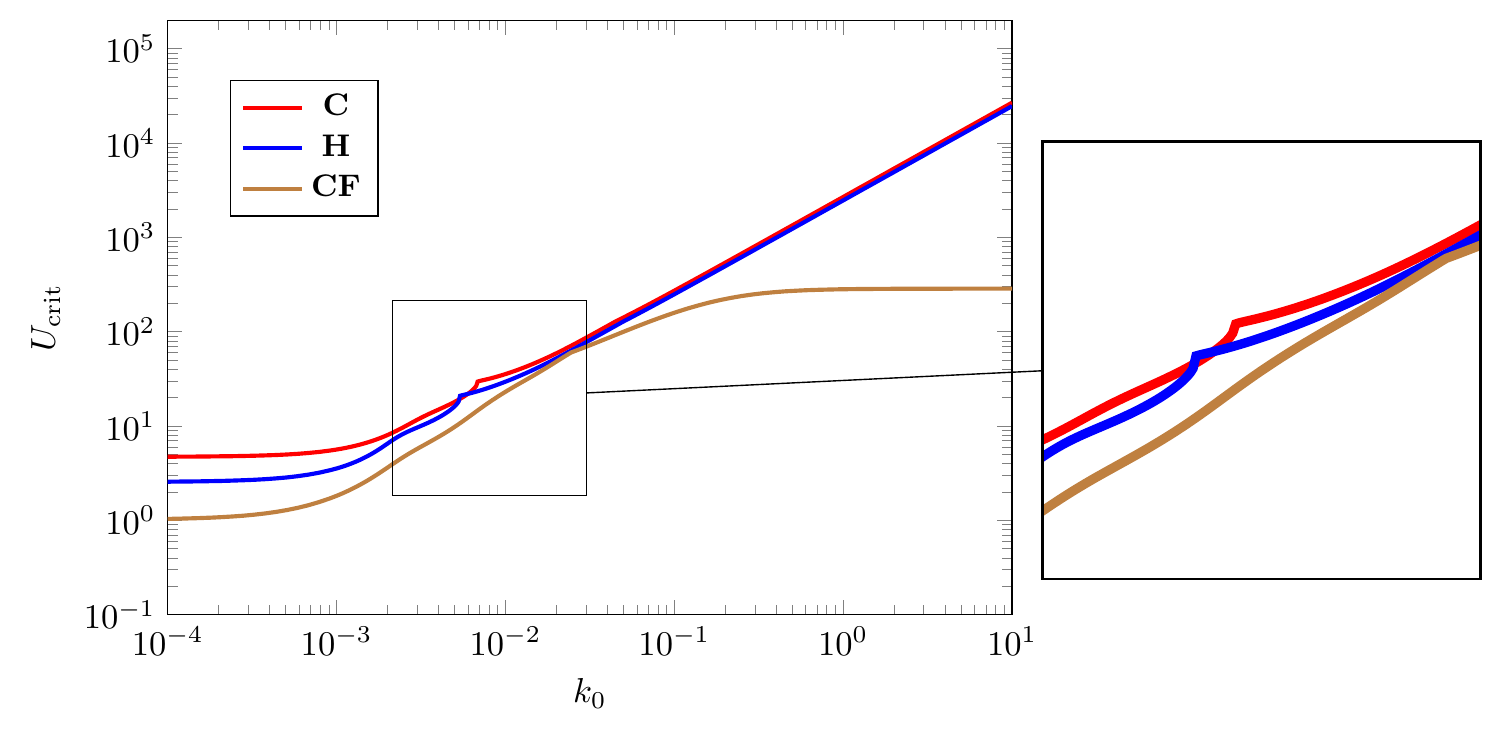}
\caption{Plot of critical $U$ values for $k_0$, $k=k_0+\beta$, $\beta=1.2\cdot10^{-4}$, $L=300$, $D=23.9$.}\label{fig3}
\end{center}
\end{figure}

Another relationship of interest is that of the damping coefficient $k=k_0+\beta$ on $\Ucrit$.
Note that in Figure \ref{fig3} the region below a given curve represents stable combinations of $k$ and $U$.  The curves themselves represent the critical value $\Ucrit$ for any given damping coefficient $k_0$.  Conversely, given a $U$ value, one can extract the necessary strength of damping to yield stability.  Note that in cases where the boundary points are completely restricted ({\bf C} and {\bf H}), the relationship between $\Ucrit$ and $k_0$ is asymptotically linear, while for the {\bf CF} configuration it is sublinear. The spurious behavior in the zoomed region on the right does not seem to be a numerical artifact, as the computations were done on a very fine scale.

\begin{figure}[H]
\begin{center}
\includegraphics[width=6in]{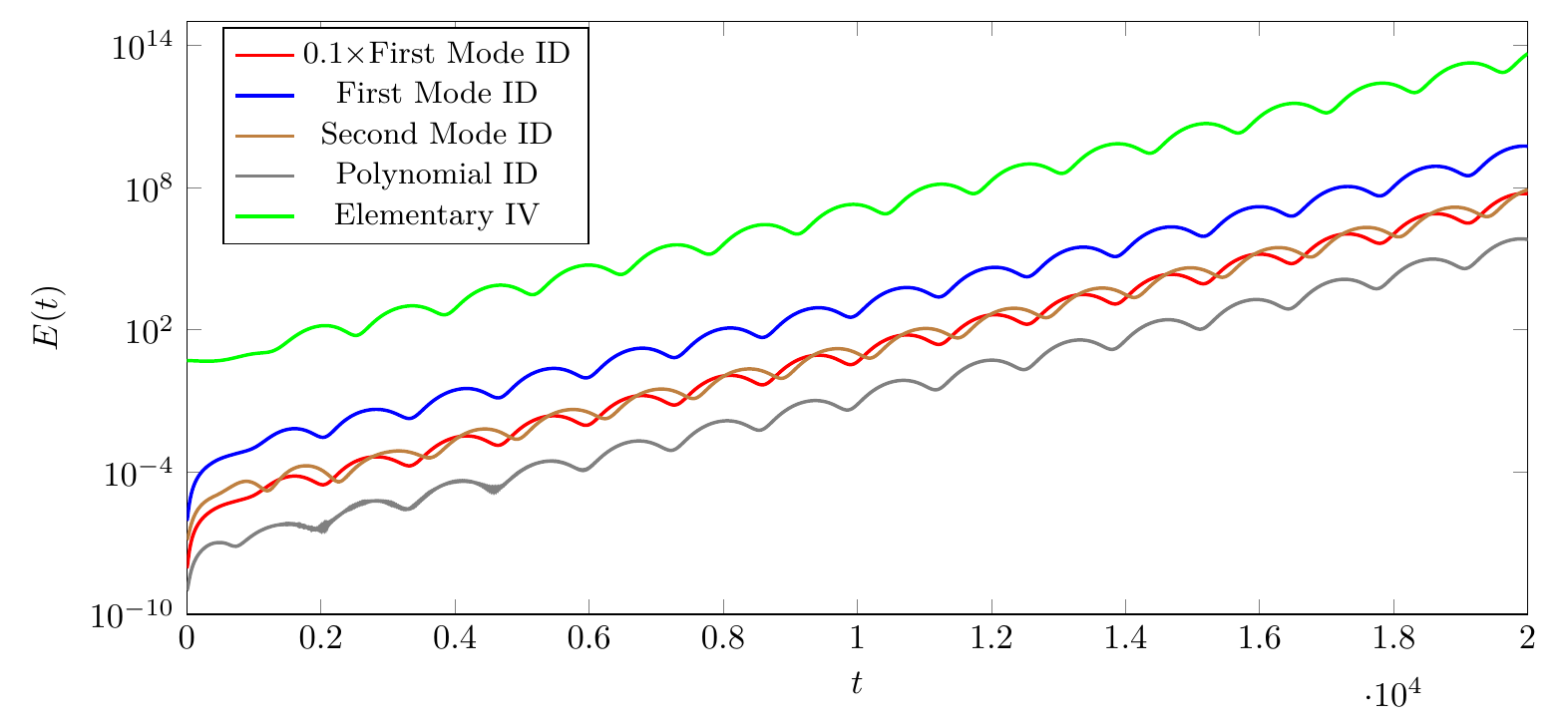}
\caption{({\bf H}) Plot of $E(t)$  with $U=5$, $k_0=0$, $\beta=1.2\cdot10^{-4}$, $L=300$, and $D=23.9$.}\label{fig4}
\end{center}
\end{figure}

We conclude this section by empirically verifying the global nature of instability across initial conditions (for a given configuration). In Figure \ref{fig4} we plot the dynamic energy as a function of time for the {\bf H} configuration with $U>U_{\text{crit}}$ (i.e., $U=5$ for $L=300$, $D=23.9$, $\beta=1.2\cdot 10^{-4}$, and $k_0=0$).  The initial conditions are taken from among the first and second {\bf H} modes (see Table \ref{tab:t1}), the polynomial initial displacement $w_0(x)=\hat x^3(1-\hat x)^3$, and the elementary initial velocity $w_1(x)=\hat x (1-\hat x)$.
We note that, in the above, all dynamic simulations reflect that $U>\Ucrit$. Each simulation, in the absence of the nonlinear restoring force ($b_2=0$), shows the dynamics eventually growing exponentially. Moreover, it is clear from the log-scale plot---since all slope profiles are the same---that each simulation is growing in time by the same exponential factor, namely, the perturbed unstable eigenvalue. Also note, as we will observe in the non-linear case ($b_2>0$), there is a transient regime at the outset of the dynamics; however, with damping provided by the flow, $\beta w_t$, the effect of the initial condition decays away after a short amount of time.

\section{Nonlinear Computational Analysis: Qualitative Properties}\label{postflutter}
In this section, we present a variety of simulations across configurations and parameter values to illustrate various aspects of the qualitative {\em post-flutter} behavior. In the simulations below, for the most part, we take the nonlinear parameter to be positive, $b_2>0$. This means, in particular, that the nonlinear restoring force is active. Thus, for any simulation with $U>U_{\text{crit}}$, we will observe convergence to a limit cycle oscillation. We can then ask about the effect of any parameter (e.g., in-axis compression $b_1>0$) on the limit cycle. 

For expediency of simulations, and also because the primary purpose of this section is qualitative, we allow parameters not of interest to be set to unity. This is to say that in Section \ref{postflutter}, we take $\beta=L=D=1$. Simulations below will consider varying $b_1, b_2, k_0$ and $U$. 

For the readers convenience, we restate the system studied in this section here:
\begin{equation}\label{Bergerplate***}
\begin{cases} w_{tt}+\partial_x^4 w + k_0w_t +[b_1-b_2||w_x||^2]w_{xx} = -  [w_t+Uw_x] \\[.1cm]
w(t=0)=w_0(\hat x);~~w_t(t=0)=w_1(\hat x)\\[.2cm]
BC(w) = \begin{cases} w(0)=w_x(0)=w(1)=w_x(1)=0 &~\text{for configuration  {\bf C } } \\[.1cm]
w(0)=w_{xx}(0)=w(1)=w_{xx}(1)=0 &~\text{for configuration  {\bf H } } \\ 
w(0)=w_x(0)=w_{xx}(1)=0,~~w_{xxx}(1)+[b_1-b_2||w_x||^2]w_x(1)=0& ~\text{for configuration} ~{\bf CF} \end{cases}\end{cases}.
\end{equation}

Implementation of the boundary condition for $w_{xxx}(L)$ in the {\bf CF} configuration was computed by utilizing a 5-point stencil for $w$ centered at $x=L$.  The nonlinear terms in \eqref{Bergerplate***} involving $b_2||w_x||^2$ were handled using the $w_x$ from the previous iterate (within the same time step) for the integral computation $||w_x||=\int_0^L|w_x|^2\,dx$. 

\subsection{Post-Flutter Examples} 
\begin{figure}[H]
\begin{center}
\includegraphics[width=2.5in]{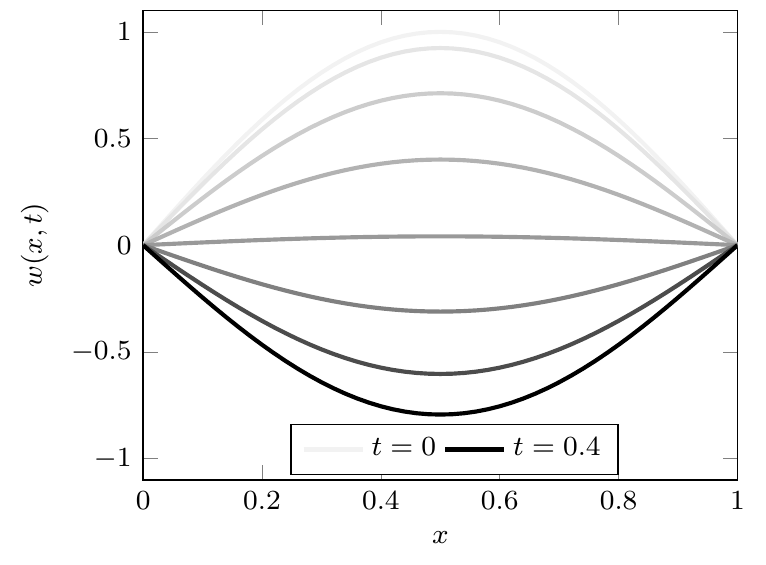}
\hspace*{0.25in}
\includegraphics[width=2.5in]{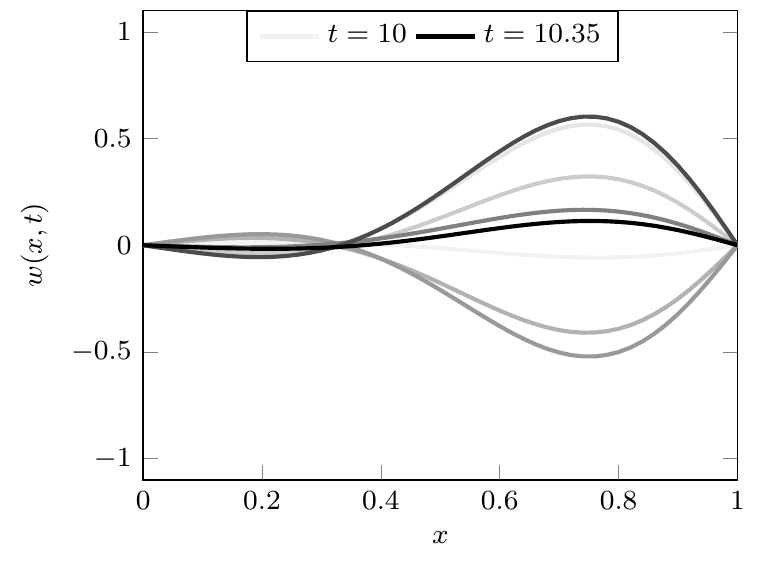}
\caption{({\bf H})  Plot of in vacuo $w(x,t)$ varying $t$, first mode as initial displacement (left); ({\bf H})  fluttering $w(x,t)$ from same initial displacement (right).}\label{figHH}
\end{center}
\end{figure}

	To ground our discussions, we begin by showing a collection of snapshots of the in vacuo \emph{linear} ($b_2=0$) beam dynamics and corresponding post-flutter dynamics.
In Figure \ref{figHH} on the left we have the in vacuo ($k_0=\beta=b_1=b_2=U=0$) dynamics with the first {\bf H} mode as an initial displacement (resulting in point-wise harmonic motion).  To the right, we have a post-flutter---$U>U_{\text{crit}}=360$ and $b_2=1$---sequence of snapshots that are demonstrative of a typical {\bf H} fluttering beam. (Note that for $U>0$, the flow is from left to right.)  

In the figures below, we provide a similar snapshot of the profile of a fluttering cantilever beam, taken with a polynomial initial displacement (see Section \ref{mainresults}), rather than a modal initial displacement. We show approximately one period of the non-transient dynamics of a beam, after it has approached a limit cycle. On the left is the natural scale, and on the right, we decrease the $y$-scale to emphasize the nodal point. For comparison, we provide time snapshots of both the first and second cantilever modes in vacuo; note the similarity between the flutter profile and a superposition of the first and second in vacuo modes. 

\begin{figure}[H]
\begin{center}
\includegraphics[width=2.5in]{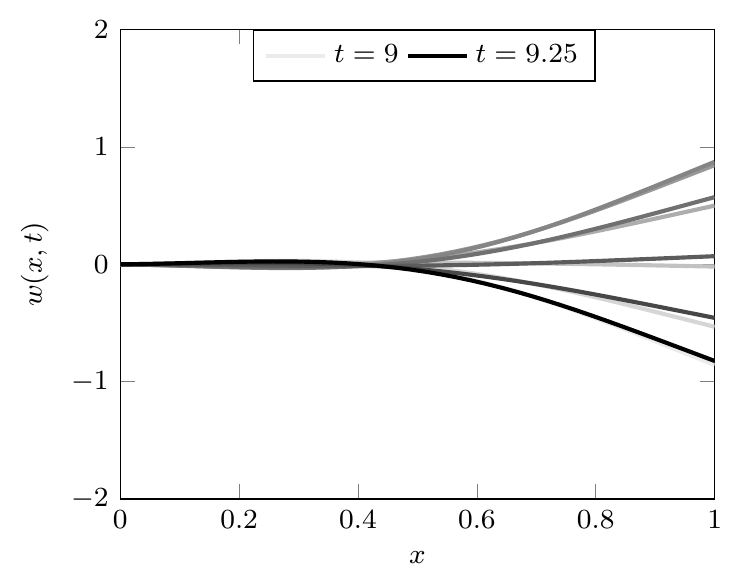}\hspace*{0.25in}
\includegraphics[width=2.5in]{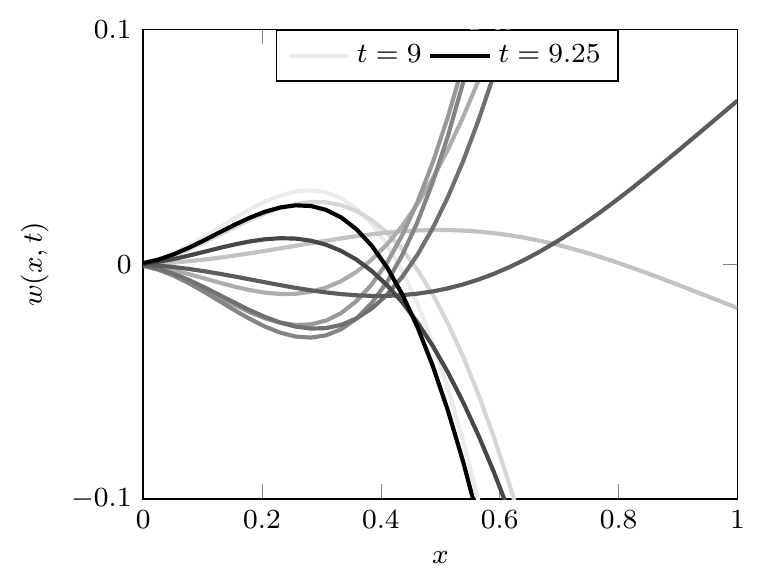}
\caption{({\bf CF}) Plot of $w(x,t)$  for $k_0=0$, $b_2=1$, $U=150$ at varying $t$ (left); to the (right) is a magnification in the transverse dimension.}\label{fig29}
\end{center}
\end{figure}

\begin{figure}[H]
\begin{center}
\includegraphics[width=2.5in]{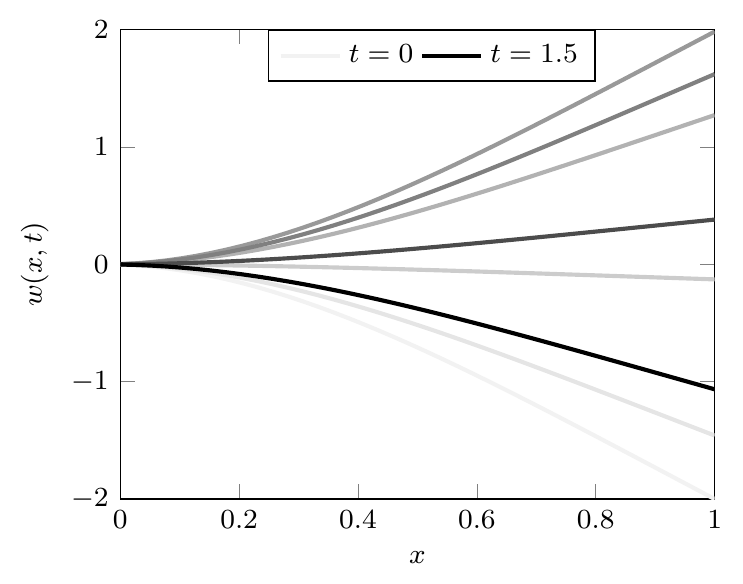}
\hspace*{0.25in}
\includegraphics[width=2.5in]{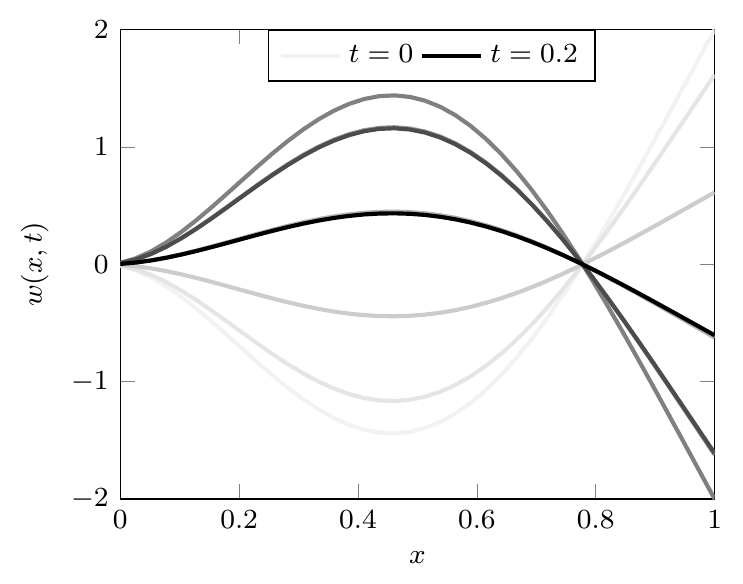}
\caption{({\bf CF})  Plot of in vacuo $w(x,t)$ at varying $t$; 1st mode as initial displacement (left), 2nd mode as initial displacement (right).}\label{fig21}
\end{center}
\end{figure}

	\subsection{Energy Plots (Linear vs. Nonlinear)} 
	
In the following section we consider the {\bf C} configuration. We  fix the initial condition and primarily consider varying $U$ around $\Ucrit$, holding other parameters fixed. We focus on the behavior of the energies to illustrate the effect of including nonlinearity in the simulation; this is to say that, when $b_2>0$, the nonlinear restoring force is active and for $U>U_{\text{crit}}$ the trajectories remain bounded (see Section \ref{long-time}) and lock into a limit cycle.

\begin{figure}[H]
\begin{center}
\includegraphics[width=5in]{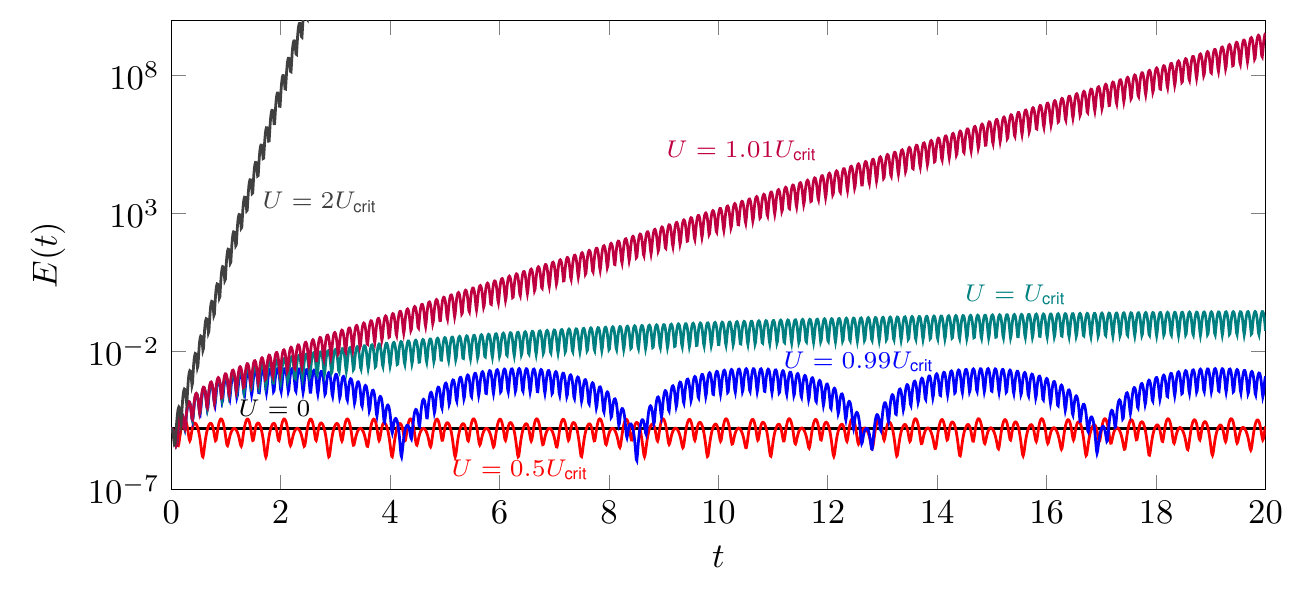}\\
\caption{({\bf C})  Plot of $E(t)$ for $k=0$, $b_1=b_2=0$, and varying $U$, $U_{\text{crit}}=135.18$.}\label{fig1l}
\end{center}
\end{figure}

First, computed energies $E(t)$ for the model \eqref{Bergerplate} with no damping present ($k_0=-\beta=-1$) and with no nonlinear or in-axis effects ($b_1=b_2=0$) are shown in Figure \ref{fig1l}.  For this choice of $k_0=-1$, we have $\Ucrit=135.18$.  Energy profiles (log-scale) for various choices of $U$ as multiples of $\Ucrit$ are given. Recall that a linear profile in the energy plot represents exponential growth or decay of the dynamics in the finite energy topology, with margin of in/stability depending on the slope of the profile.  Note the exponential growth in energies for $U>\Ucrit$ while $E(t)$ stays bounded for $0<U<\Ucrit$, and is constant for $U=0$.  

\begin{figure}[H]
\begin{center}
\includegraphics[width=5in]{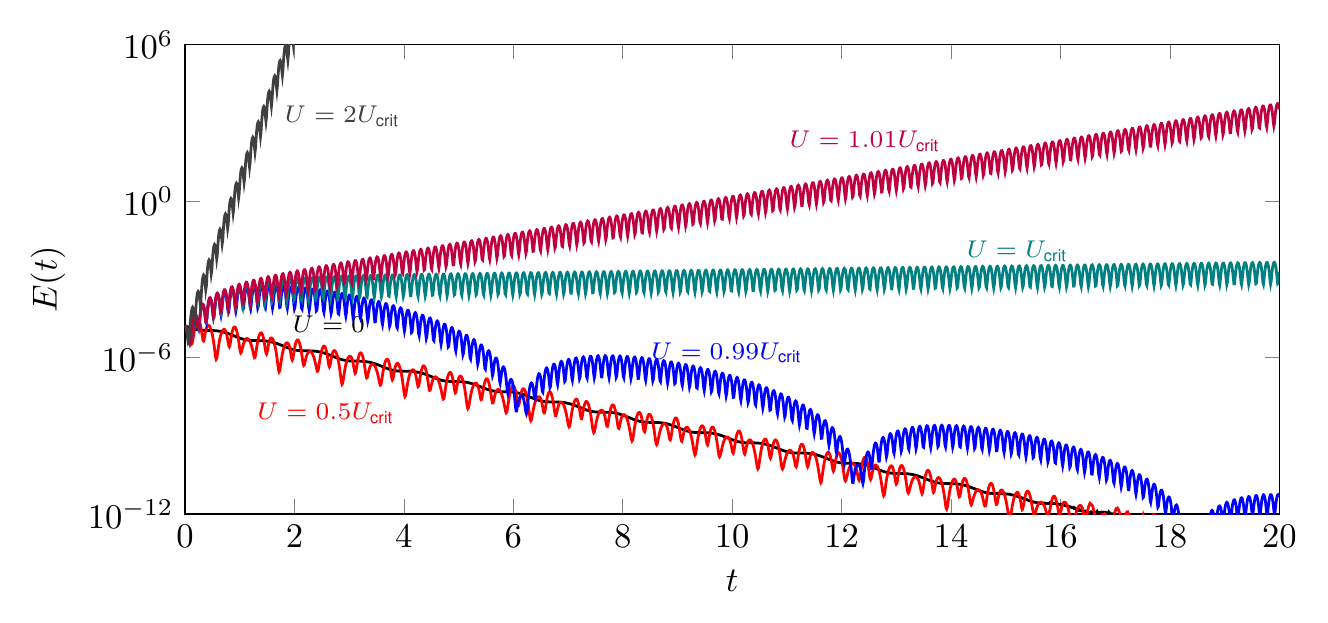}\\
\vspace*{-0.2in}
\caption{({\bf C}) Plot of $E(t)$  for $k=1$, $b_1=b_2=0$,  and varying $U$, $\Ucrit=135.9$.}
\label{fig2l}
\end{center}
\end{figure}

In the next simulation the damping parameter is set to $k_0=0$ so $k=1$, and thus there {\em is} damping present from the flow (but no structural or imposed damping) in \eqref{Bergerplate0}. This results in an approximate critical flow velocity of $\Ucrit=135.9$ (the presence of damping perturbs the critical flow velocity corresponding to the onset of instability by roughly $0.5\%$).  Computations again were performed  for various choices of $U$, with energies shown in Figure \ref{fig2l}.  Note the effect that the damping has on the decay of energies for all $U<\Ucrit$ (as $E(t)\to 0$ when $t\to\infty$)---each curve for $U<\Ucrit$ has a linear profile (or envelope) with negative slope, indicating exponential decay.  Energy values for $U\ge \Ucrit$ still grow exponentially.

	Finally, in Figure \ref{fig4l}, the clamped, nonlinear beam (with $b_1=0$ and $b_2=1$) is considered and the nonlinear energy $\mathscr{E}(t)$ is shown.  Note that $\mathscr{E}(t)$ (and hence $E(t)$, due to \eqref{nonen}) now remains bounded for all supercritical (unstable) velocities $U>\Ucrit=135.9$, demonstrating the stability that is induced by nonlinear restoring force. From the point of view of trajectories, each remains bounded for all time, with global-in-time bound dependent upon $U$. For $U<\Ucrit$, we note that the {\em nonlinearity does not dramatically affect the stability observed in the linear case}---consistent with the semilinear nature of the nonlinearity.

\begin{figure}[H]
\begin{center}
\includegraphics[width=5in]{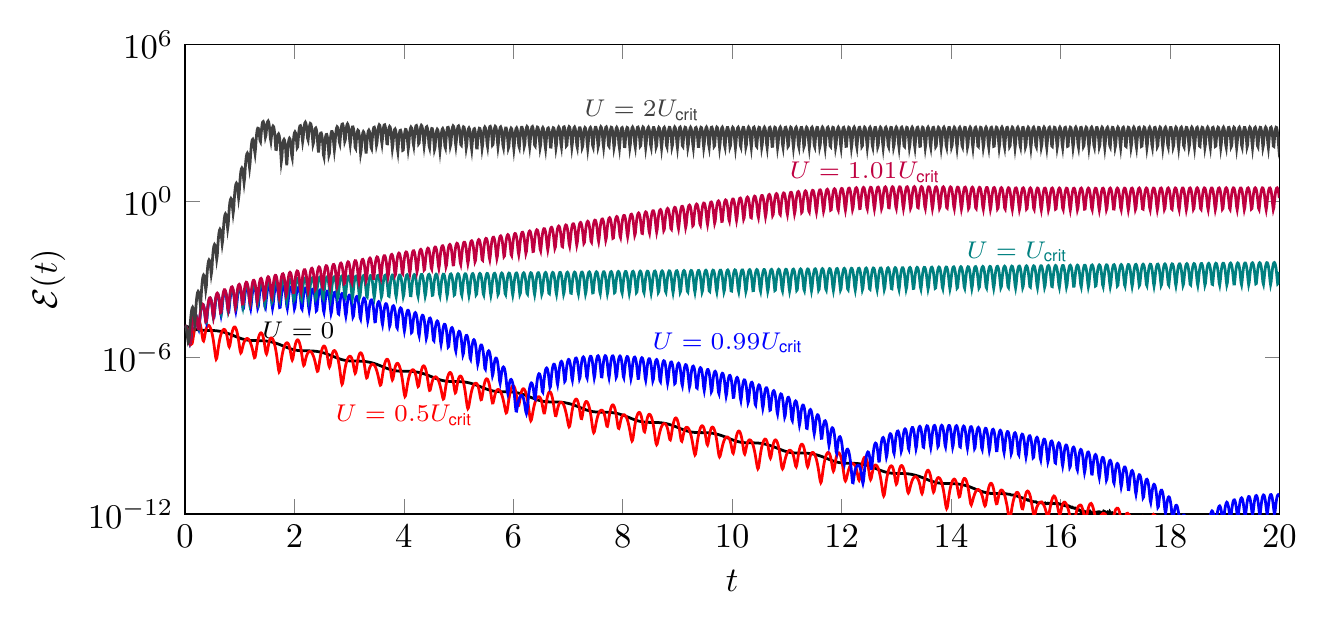}\\
\vspace*{-0.2in}
\caption{({\bf C}) Plot of $\mathscr{E}(t)$  for $k=1$, $b_1=0$, and $b_2=1$, varying $U$, $\Ucrit=135.9$.}
\label{fig3l}
\end{center}
\end{figure}

	To show the detailed difference between the energies $E(t)$ and $\mathcal{E}(t)$, Figure \ref{fig5l} shows both quantities for $k_0=0$, $b_1=0$, and $b_2=1$, with $U=2\Ucrit$.  
\begin{figure}[H]
\begin{center}
\includegraphics[width=5in]{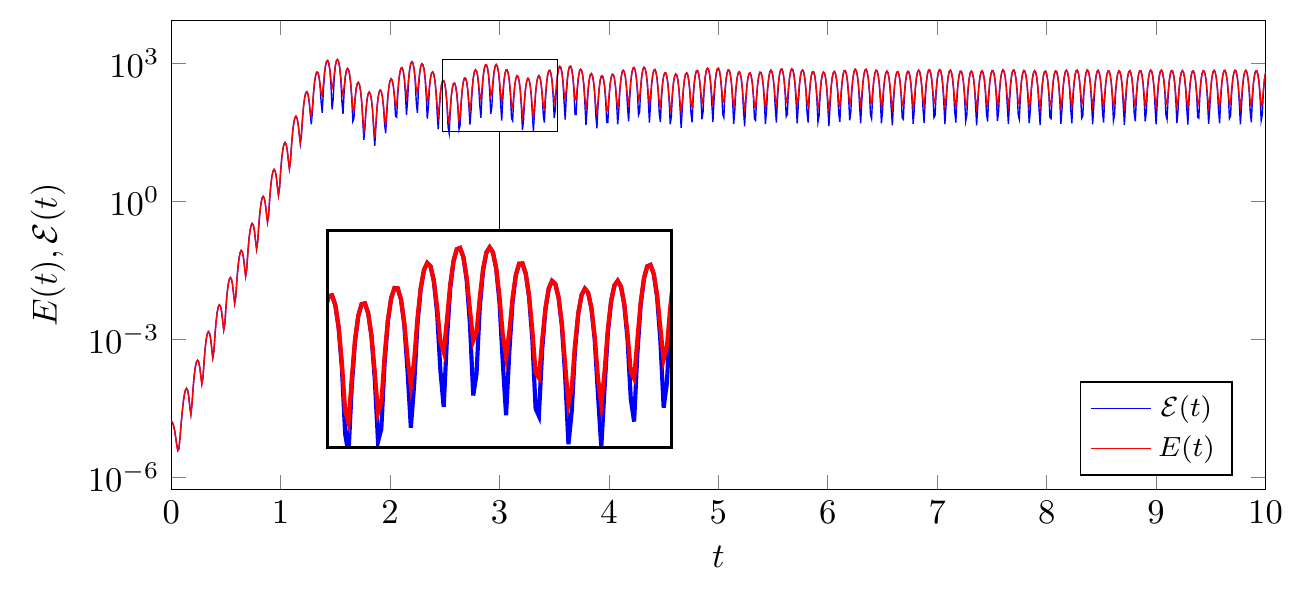}\\
\vspace*{-0.2in}
\caption{({\bf C}) Plot of energies  for $k=1$, $b_1=0$, $b_2=1$, and $U = 2\Ucrit$.}
\label{fig4l}
\end{center}
\end{figure}

\subsection{Blow Up for Nonlinear {\bf CF} with Linear Boundary Conditions}\label{badbcblowup}
We recall that, for conservation of energy to hold (or, equivalently, to derive the equations of motion through Hamilton's principle) in the {\bf CF} configuration, the boundary conditions \eqref{badBC}  must be altered; in the panel configurations {\bf C} and {\bf H}, standard linear boundary conditions \eqref{Bergerplate} are correct. (For more discussion, see \cite{HTW}.)

A natural question, then, for the nonlinear cantilever ($b_2>0$) is to ask about the stability of the \eqref{Bergerplate0}  model in the presence of {\em standard linear} clamped-free boundary conditions (as in \eqref{Bergerplate***}). We emphasize that this combination is non-physical, as there is no {\em conservation of energy} associated to the dynamics. Below, we achieve arbitrary growth (in time) of displacements (or energies) for finite difference solutions to:
\begin{equation}
\begin{cases} w_{tt}+\partial_x^4w-\|w_x\|^2w_{xx}  = 0 \\
w(t=0)=0;~~w_t(t=0)=cx \\
w(0)=w_x(0)=0;~~w_{xx}(1)=0,~~\partial_x^3w(1)=0, \end{cases}
\end{equation}
for $c$ sufficiently large.

Specifically, if $c=12$, the dynamics exhibit periodic behavior; when $c=13$, the dynamics grow with exponential rate. When the boundary conditions are adjusted appropriately (i.e., $\partial_x^3w(1)=||w_x||^2w_x(1)$), the nonlinear energy $\mathscr E(t)$ is nearly {\em perfectly conserved}.

\begin{figure}[H]
\begin{center}
\includegraphics[width=5in]{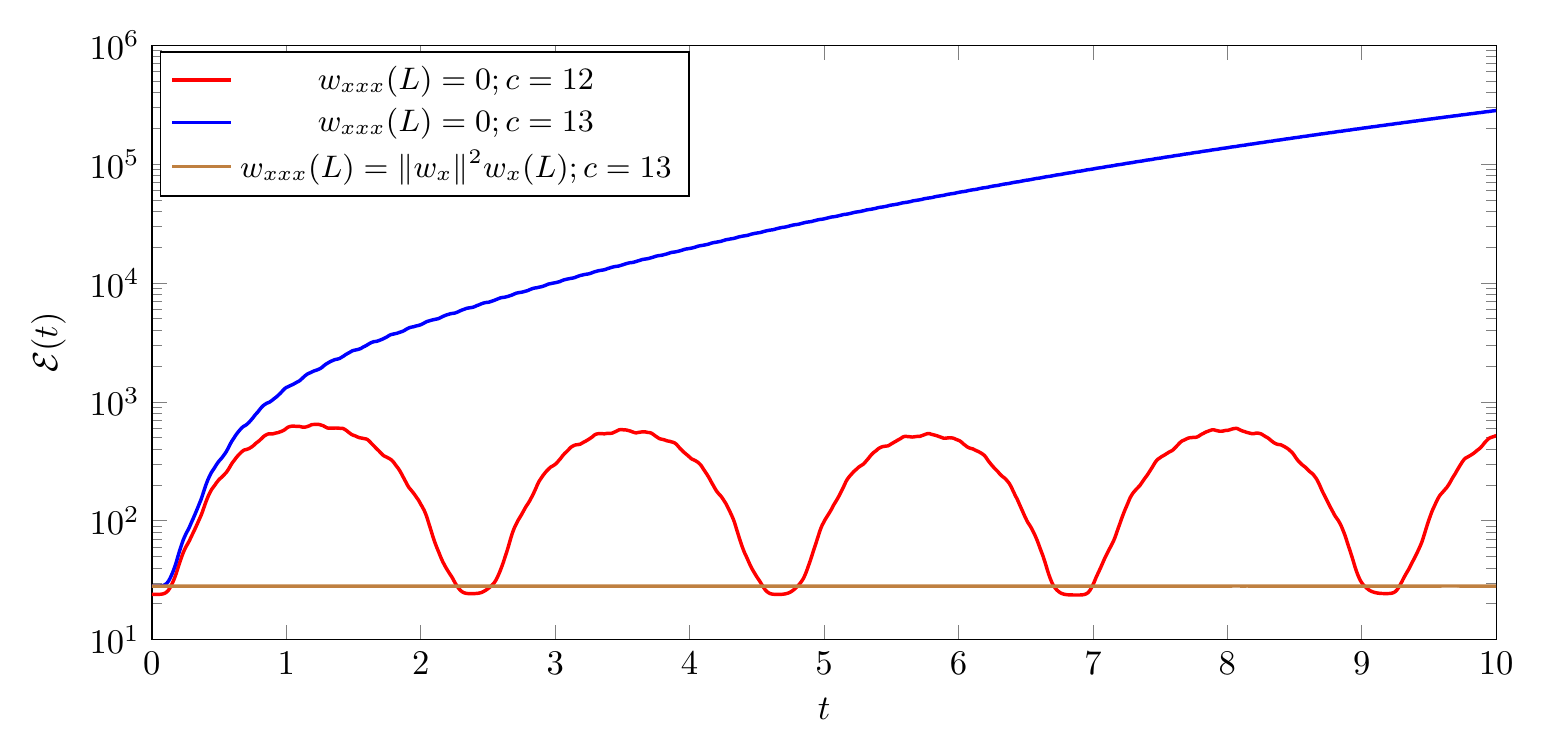}\\
\vspace*{-0.2in}
\caption{({\bf CF}) Plot of nonlinear energies for comparing linear and nonlinear boundary conditions.}
\label{fig5l}
\end{center}
\end{figure}

	\subsection{Convergence to Limit Cycles}\label{LCOs}
	In this section we demonstrate some properties of trajectories that converge to a limit cycle oscillation, as described above. 
	
\begin{figure}[H]
\begin{center}
\includegraphics[width=5in]{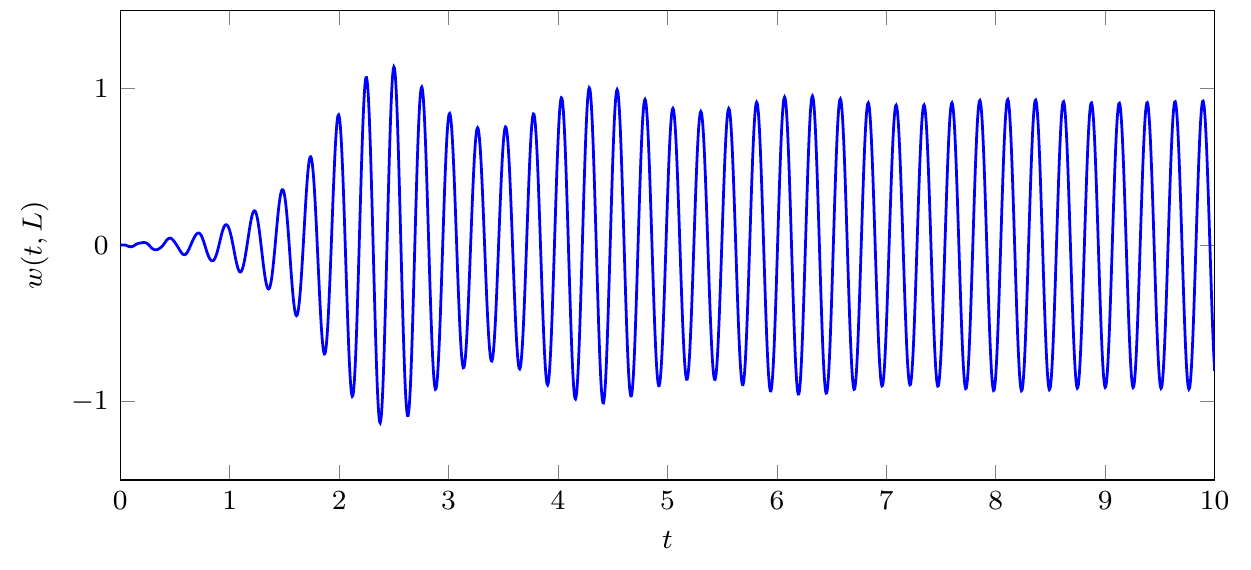}\\
\vspace*{-0.2in}
\caption{({\bf CF})  Plot of $w(t,L)$  for  $k = 1$, $b_1=0$, $b_2=1$, and $U = 2\Ucrit$.}
\label{fig6l}
\end{center}
\end{figure}

Figure \ref{fig6l} shows the endpoint displacement a cantilevered beam {\bf CF} with supercritical parameter values.  It is clear that the envelope of the oscillations grows significantly and until around $t=1.5$, and then undergoes an exponential decay to a nontrivial limit cycle behavior.

\begin{figure}[H]
\begin{center}
\includegraphics[width=5in]{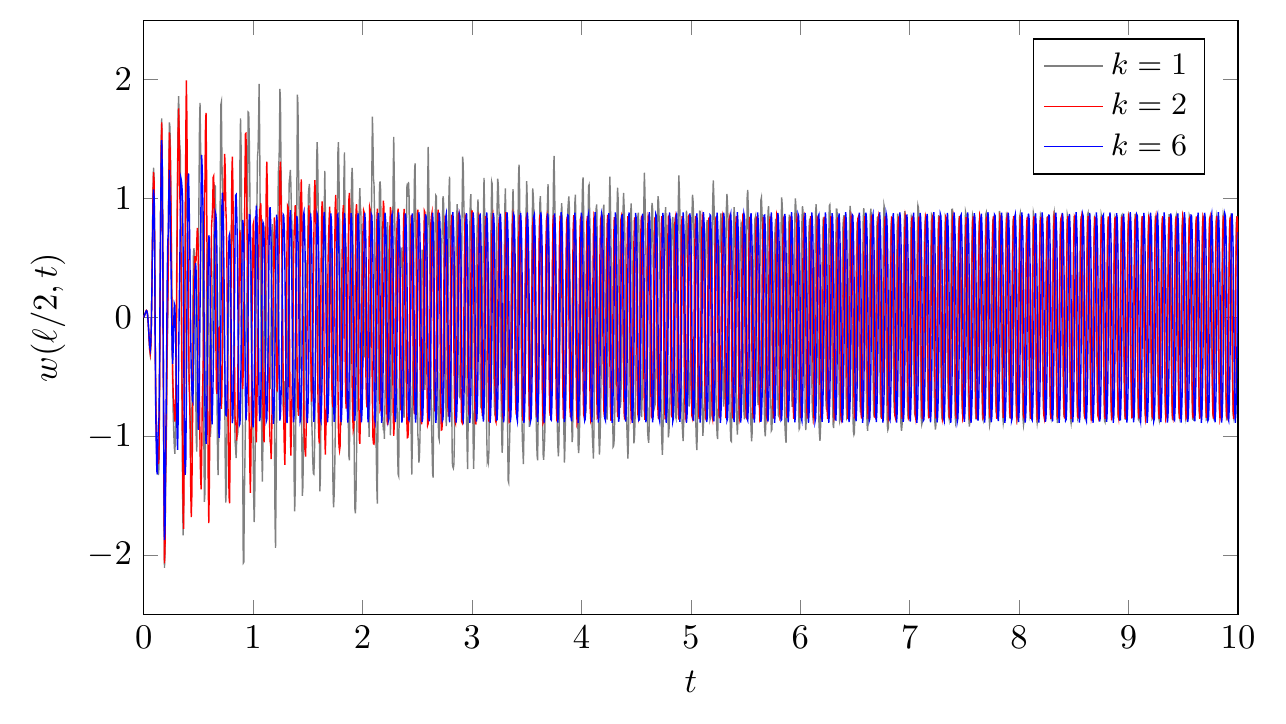}\\
\vspace*{-0.2in}
\caption{({\bf C})  Plot of $u$ at beam midpoint  for $U=5000$, $b_1=20$, $b_2=1$, varying $k$.}
\label{fig11l}
\end{center}
\end{figure}

In Figure \ref{fig11l}, plots of the beam midpoint displacement are shown for a very large flow velocity ($U=5000$), significant in-axis compression parameter ($b=50$), and $b_0=1$ for selected values of the damping parameter $k=k_0+1$.  The sensitivity of the dynamics to the damping parameter can be seen by noting the relative rate of initial decay of energy as $k$ increases---note the quick decay to the limit cycle for larger values of $k$, though the limit cycle itself is preserved across damping values.   

	\subsection{Uniformity of Limit Cycle}

\begin{figure}[H]
\begin{center}
\includegraphics[width=5in]{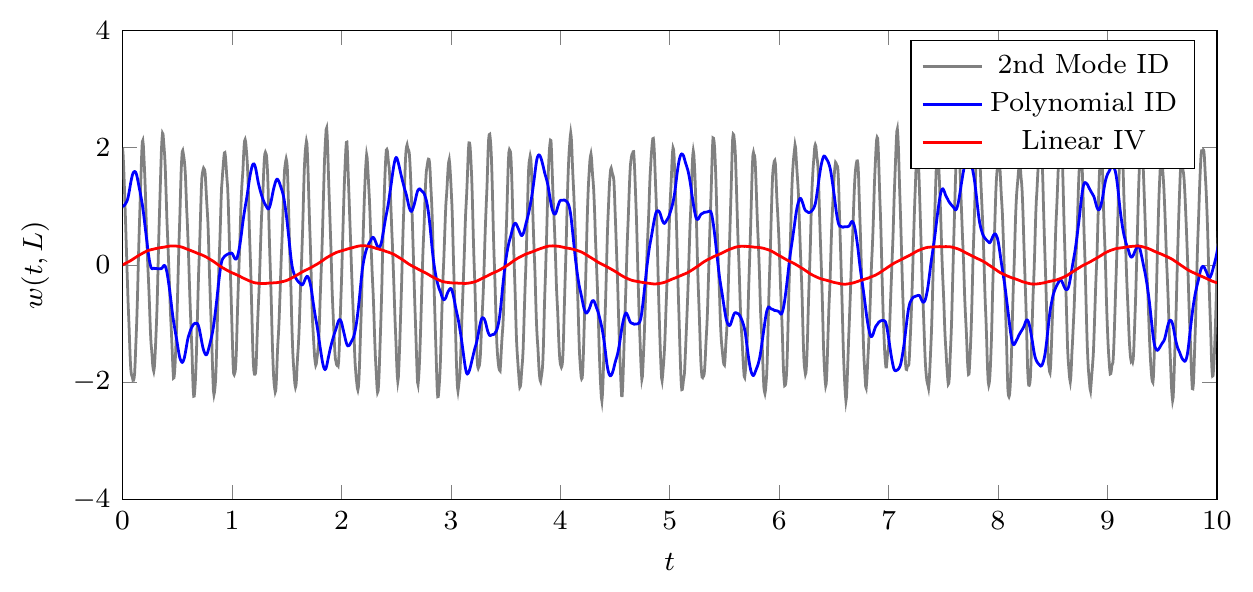}\\
\caption{({\bf CF}) Plot of in vacuo tip displacement $w(t,L)$; varying initial configuration.}\label{fig26}
\end{center}
\end{figure}

\begin{figure}[H]
\begin{center}
\includegraphics[width=5in]{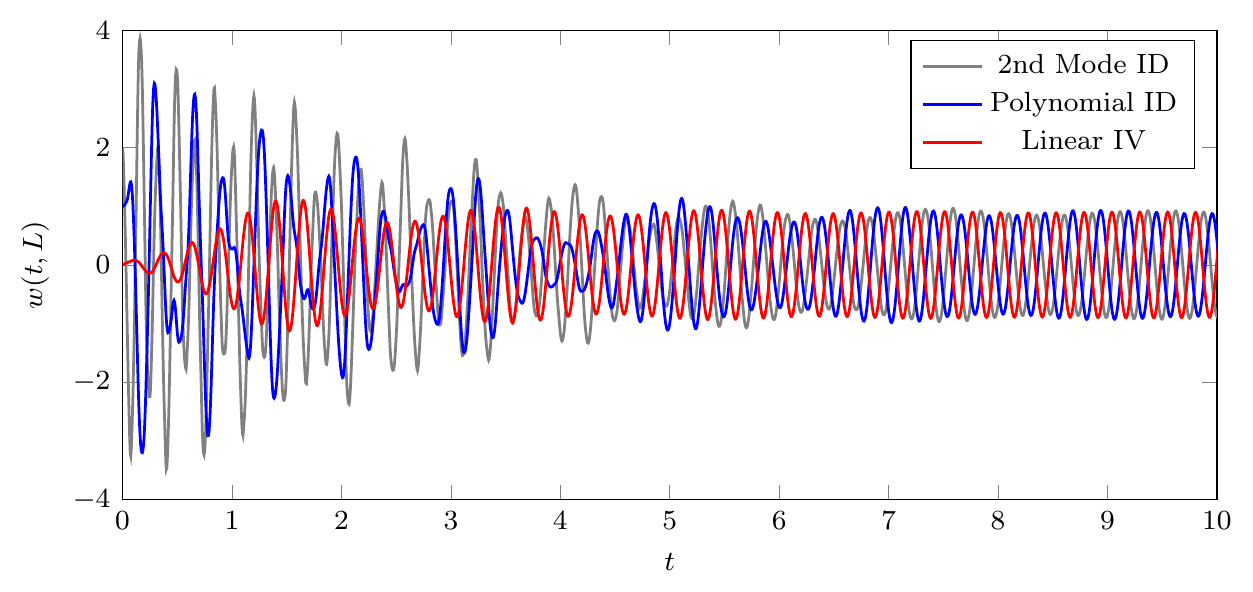}\\
\caption{({\bf CF}) Plot of $w(t,L)$ for $k_0=0$, $b_2=1$, $U=150$; varying initial configuration.}\label{fig23}
\end{center}
\end{figure}

	Next, we consider the {\bf CF} configuration  and provide tip profile at $x=L=1$ for three different initial configurations:
\begin{itemize}
  \item \textbf{[2nd Mode ID]}
$w(0,x)=s_2(x) =  [\cos(\kappa_2x)-\cosh(\kappa_2x)]-\mathcal C_2[\sin(\kappa_2x)-\sinh(\kappa_2x)], ~w_t(0,x)=0$,  \vskip.05cm
where $\kappa_2\approx 4.6941$ is the second Euler-Bernoulli cantilevered mode number (with $L=1$) and
 
 $\mathcal C_2 = \left[\dfrac{\cos(\kappa_2)+\cosh(\kappa_2)}{\sin(\kappa_2)+\sinh(\kappa_2)}\right]\approx 1.0185$; 
 \item \textbf{[Polynomial ID]}\quad
$w(0,x)=-4x^5+15x^4-20x^3+10x^2$, ~$w_t(0,x)=0$; 
\item \textbf{[Linear IV]}\quad $w(0,x)=0$, $~w_t(0,x)=x$. 
\end{itemize} 
First, Figure \ref{fig26} represents the nonlinear ($b_2=1$), in vacuo ($\beta=k_0=0$) tip displacements with the three initial configurations described above.
Figure \ref{fig23} represents fluttering dynamics across the three different initial configurations. 

In these cases, we observe convergence to the ``same" limit cycle for each of the initial conditions considered. Finally,  note that, although we seem to see convergence to the same LCO, the {\em transient regime} and ``time to convergence" are certainly affected by the choice of initial configuration. 

	\subsection{Convergence to Nontrivial Equilibrium} 
	
\begin{figure}[H]
\begin{center}
\includegraphics[width=5in]{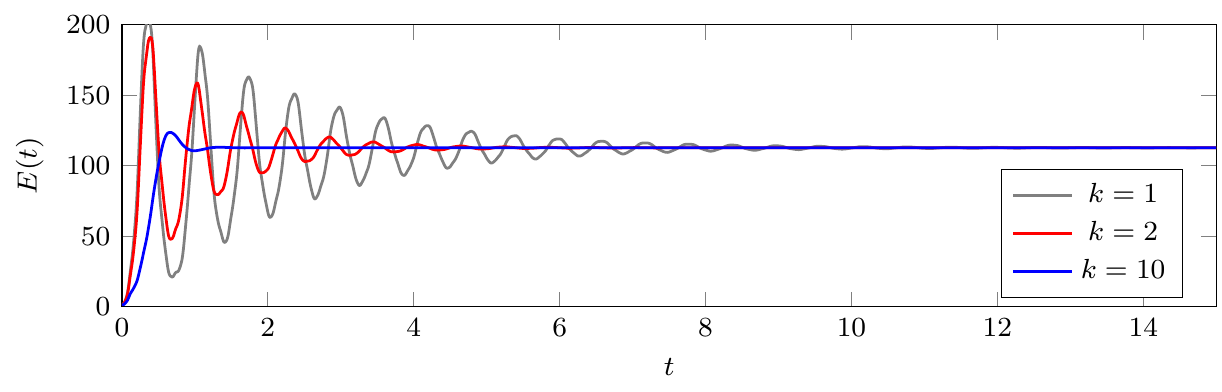}\\
\vspace*{-0.2in}
\caption{({\bf C}) Plot of $E(t)$  for  $U=100$, $b_1=50$, $b_2=1$, and varying $k$.}
\label{fig8l}
\end{center}
\end{figure}

\begin{figure}[H]
\begin{center}
\includegraphics[width=5in]{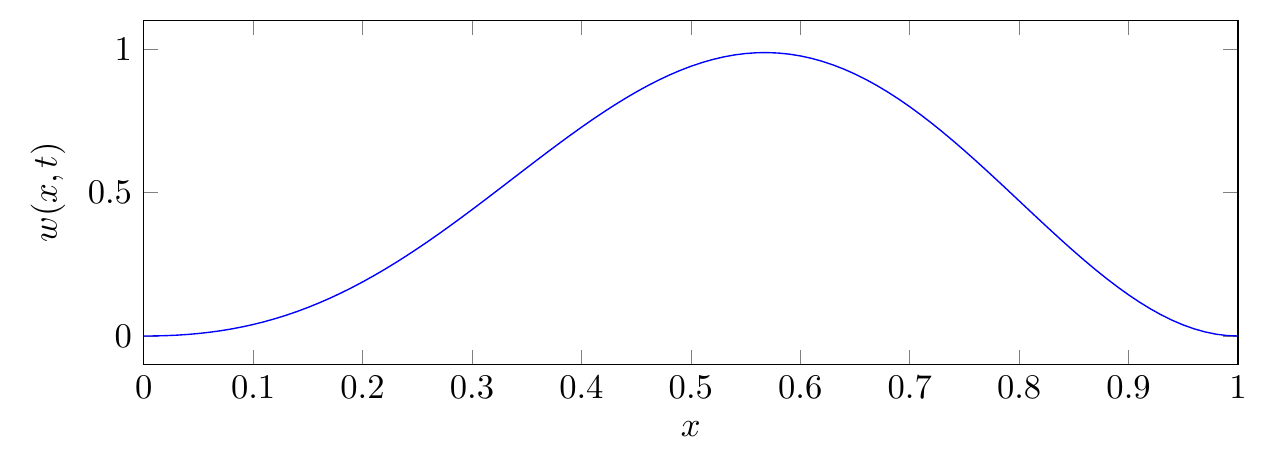}\\
\vspace*{-0.2in}
\caption{({\bf C}) Plot of steady-state beam displacement  for $U=100$, $b_1=50$, $b_2=1$, and varying $k$.}
\label{fig9l}
\end{center}
\end{figure}

	For certain parameter combinations, the presence of excess in-plane compression leads to a ``buckled'' plate/beam configuration (a bifurcation, for the static problem).  For $b_2=1$, the parameter $b_1=50$ is large enough for  $U=100$ to impart this behavior.  The transient behavior decays more rapidly to the nontrivial steady state for larger values of $k$.  A plot of the energy $E(t)$ for different values of $k$ is given in Figure \ref{fig8l}.  Regardless of the choice of $k$, all simulations for $U=100$, $b_2=1$, and $b_1=50$ converge to the same nontrivial steady state and same linear energy level set.  A plot of the nontrivial steady state $w$ (the buckled beam displacement) is given in Figure \ref{fig9l}.

\begin{figure}[H]
\begin{center}
\includegraphics[width=5in]{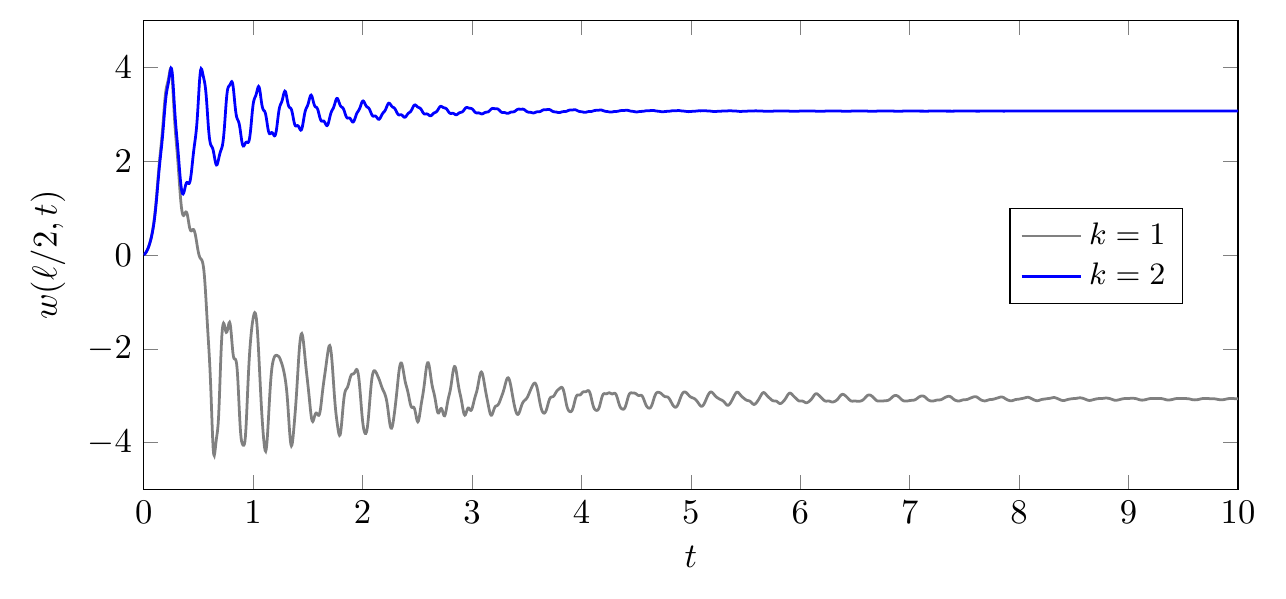}\\
\vspace*{-0.2in}
\caption{({\bf C})  Plot of $w$ at beam midpoint  for $U=100$, $b_1=100$, $b_2=1$, and varying $k$.}
\label{fig10l}
\end{center}
\end{figure}

It is also possible to observe different choices of damping $k$ resulting in convergence to different steady states.  For $U=100$, $b_2=1$, and $b_1=100$, the choices $k=1$ and $k=2$ actually produce nontrivial steady states that are negatives of one another.  In Figure \ref{fig10l} the midpoint displacement is plotted for these two cases.  
	\subsection{Non-simple Limit Cycle}

\begin{figure}[H]
\begin{center}
\includegraphics[width=5in]{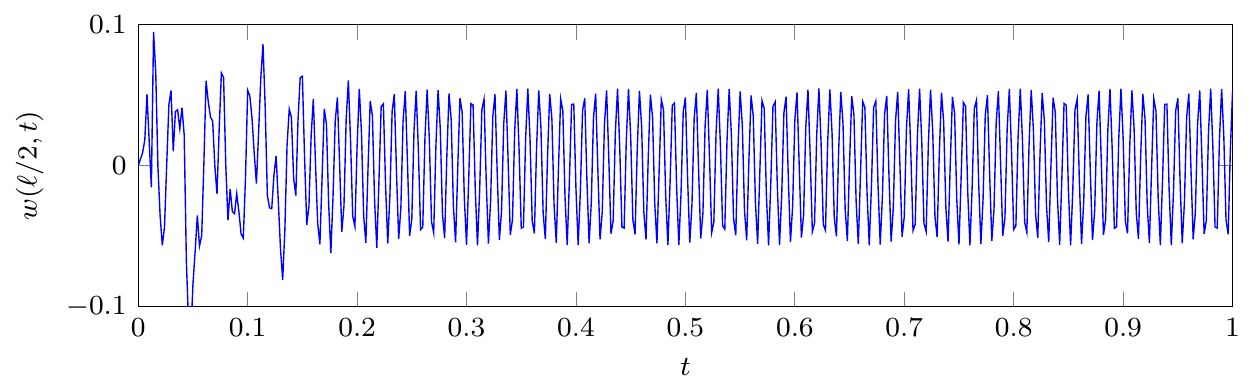}\\
\vspace*{-0.2in}
\caption{ ({\bf C}) Plot of $w$ at beam midpoint  for $U=5000$, $b_1=5000$, $b_2=5000$, $k=101$.}
\label{fig12l}
\end{center}
\end{figure}

It is possible to induce several different phenomena by manipulating parameter values for the nonlinear dynamics.  In Figure \ref{fig12l} the midpoint displacement is shown for $U=5000$, $k=100$, $b_1=5000$, and $b_2=1000$.  Note the initial transient dynamics are damped out quickly and the dynamics converges to a non-simple limit cycle. 

	\subsection{Chaos with In-plane Tension}\label{chaosbaby}
	
Here we consider the following initial configurations for the nonlinear {\bf H} beam taken with substantial in-plane compression. 
$$w(t=0,x)=w_0(x)=\epsilon \sin(2\pi \hat x),~~~~w_t(0,x)=w_1(x)=\hat x (1-\hat x); ~~~\epsilon \ge 0.$$
We note from \cite{dowellchaos} that this setup is a candidate to produce a deterministic, chaotic dynamical system. (In the reference \cite{dowellchaos}, the system \eqref{Bergerplate0} is reduced down to a two modal system, and the qualitative properties---convergence to equilibria, convergence to a simple limit cycle, convergence to a non-simple limit cycle, and ``chaos"---are mapped in the $(b_1,U)$ parameter space.)

We produce an interesting example below with large $U$. We note that for the given value of $b_1$ below, we are past the {\em Euler buckling} bifurcation point. 

	\begin{figure}[H]
\begin{center}
\includegraphics[width=5in]{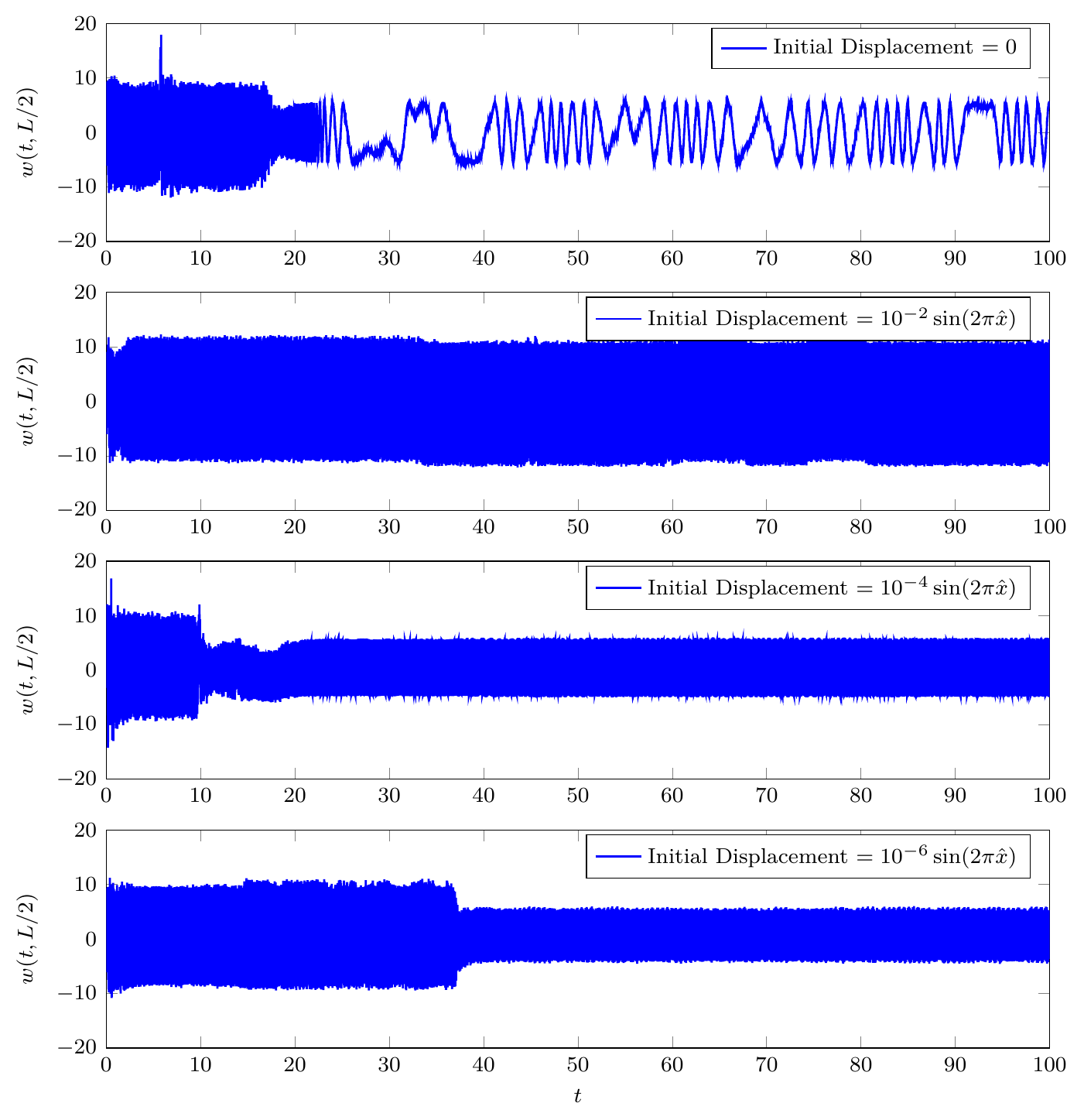}\\
\vspace*{-0.2in}
\caption{({\bf H}) Plot of $w$ at beam midpoint  for $U=200$, $b_1=2000$, $b_2=1$, $k_0=1$, and $\beta=1$.}
\label{figll}
\end{center}
\end{figure} 

First, from the four profiles corresponding to $\epsilon \searrow 0$, we see no discernible notion of convergence of the dynamics. Secondly, we see no periodic behavior corresponding to the initial configuration with $\epsilon=0$. Additionally, it is not clear that the other configurations involving $\epsilon>0$ converge to a true limit cycle. In these cases, it is clear there is no uniformity of end behavior across initial configurations, as was the case in the previous section.

\begin{remark} 
We continued simulations for smaller values of $\epsilon$, and it was clear that  no trend in magnitude or period of the oscillations emerged, even at very fine scales and long run-times.\end{remark}

\begin{figure}[H]
\begin{center}
\includegraphics[width=5in]{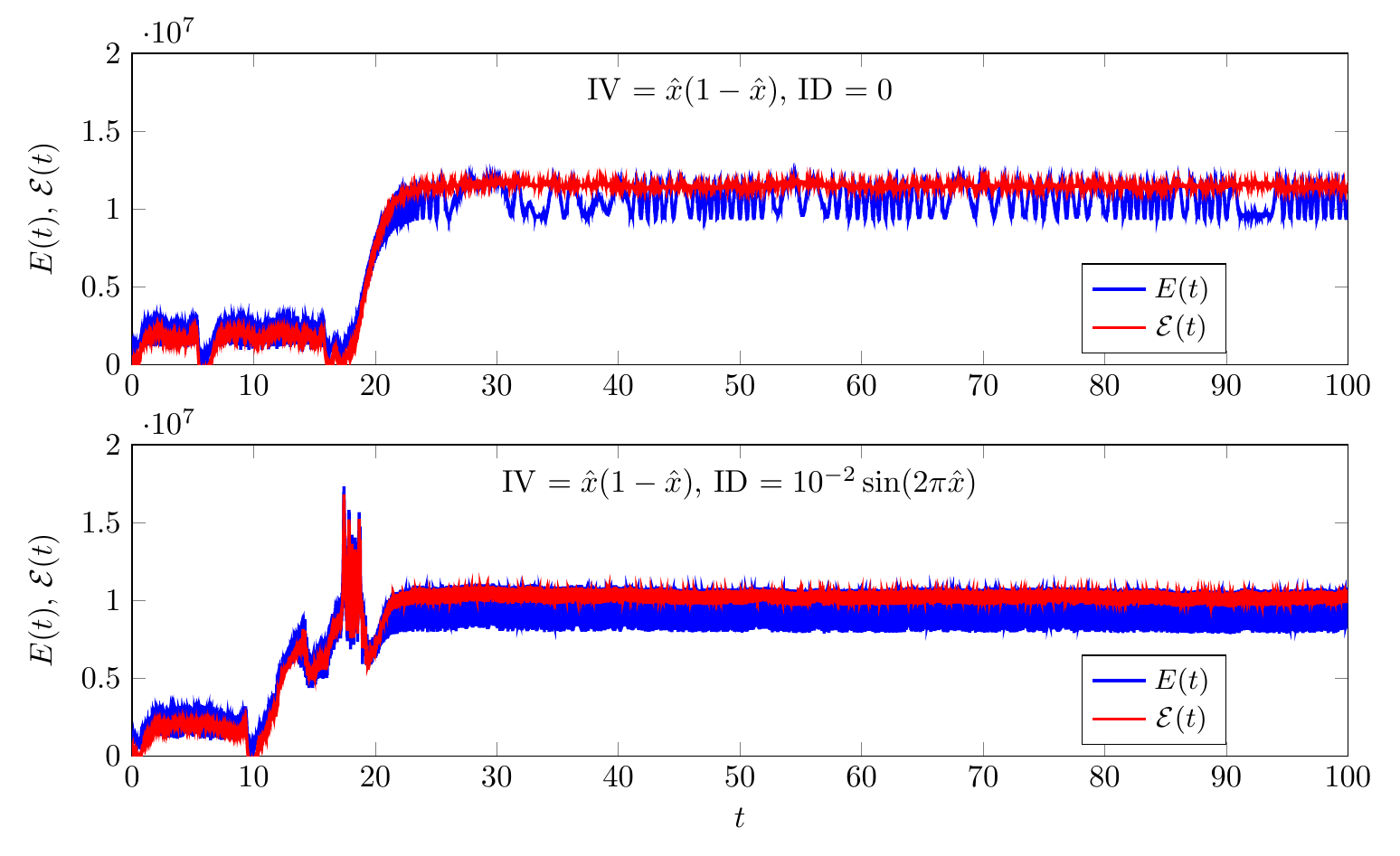}\\
\vspace*{-0.2in}
\caption{({\bf H})  Plot of $E(t), \mathcal{E}(t)$  for $U=200$, $b_1=2000$, $b_2=1$, $k_0=1$, and $\beta=1$.}
\label{figlll}
\end{center}
\end{figure}

The only consistent feature across the above simulations corresponds to their energy curves. All energy plots show the same features present in Figure \ref{figlll}; namely a transient period followed by an energy ``bump" up to a higher level. The transition is smooth. 

		\subsection{Effect of Magnitude of $b_2$}

\begin{figure}[H]
\begin{center}
\includegraphics[width=5in]{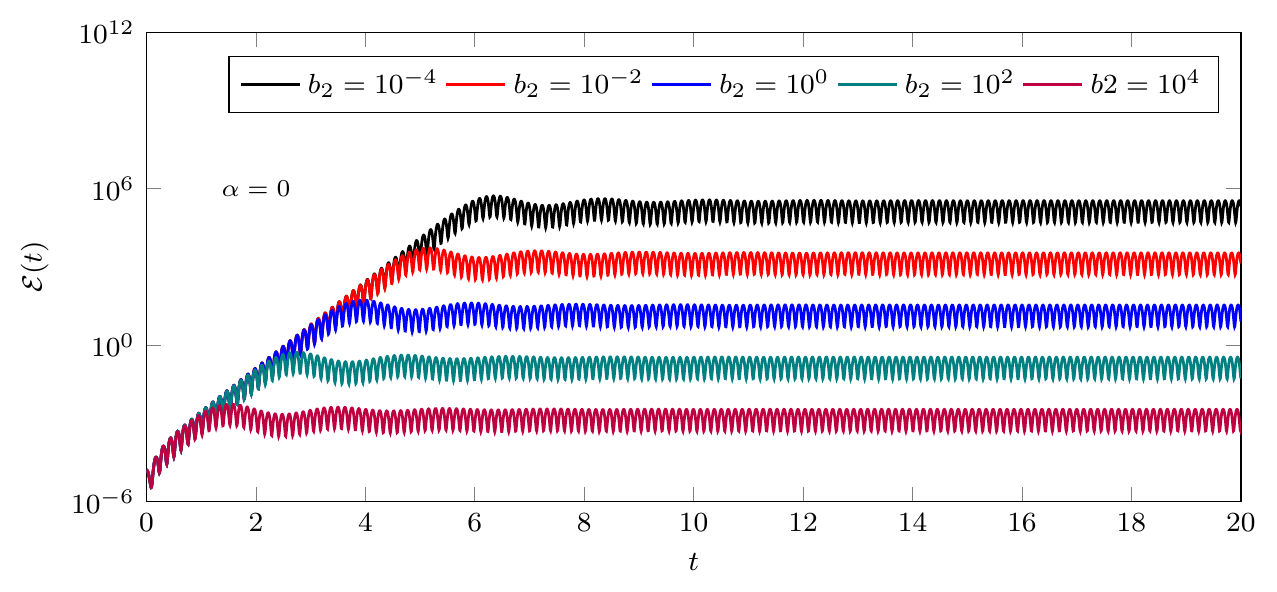}\\
\vspace*{-0.2in}
\caption{({\bf C}) Plot of $\mathscr{E}(t)$  for $k_0=0$, and $U=150$, varying $b_2$.}
\label{fig11}
\end{center}
\end{figure}

The effect of increasing the nonlinear parameter $b_2$ was studied at a supercritical flow velocity ($U >\Ucrit$) for beam.  In Figure \ref{fig11}, energy profiles for several different choices of $b_2$ are given for the parameter $k_0=0$  and $U=150$ ($\Ucrit=135.97$).  Note that  as $b_2$ increases, the energy plateau decreases (for an initial configuration fixed across all simulations).

\section{Synopsis of Main Numerical Observations}\label{end}
Here we briefly outline the main observations in Sections \ref{linear} and \ref{postflutter} above, and succinctly connect to the theory in Section \ref{sec:theory}. 

\vskip.1cm
\noindent{\bf Modal Approximation Works Well}
\vskip.1cm
In Section \ref{modal} we describe the ``modal" approach used commonly in engineering to extract the {\em flutter point} (or point of instability, with respect to critical parameters). We demonstrate in Section \ref{predict} that---for physical parameters---the modal approach faithfully replicates the onset of instability with various parameters obtained empirically (through direct numerical simulations of the PDE dynamics). This approximation is very good across configurations, and across large ranges of relevant parameters. 

\vskip.1cm
\noindent{\bf Extensible Nonlinearity Controls Non-conservative Effects}
\vskip.1cm
In various simulations above we have confirmed what we established theoretically: the extensible nonlinearity featured here in \eqref{Bergerplate} ($b_2>0$) is sufficient to control the non-dissipative (non-conservative) effects introduced by the flow through piston theory \eqref{linpist}. This is independent of configuration, of course making the necessary change to the boundary condition in the {\bf CF} configuration---see Section \ref{badbcblowup}. When the damping due to piston-theory is taken into account, the overall bound on trajectories, as the transient dynamics decay, seems independent of the initial data (i.e., its size and type), and the magnitude of the non-transient behavior (typically an LCO) depends on the intrinsic parameters, such as $b_2,U,\beta,L$, etc. 

\vskip.1cm
\noindent{\bf Onset of Instability is {\em Linear}}
\vskip.1cm
Through our simulations, we have observed that below the point of onset (for instance, when $U<U_{\text{crit}}$ in a particular configuation, the dynamics are effectively not altered by the presence of the (higher order) semilinear effect of extensibility. This means that the emergence of instability (bifurcation) is not really effected by the value of $b_2$.  In this way, we may claim that the onset of flutter is a purely linear phenomenon, regarding the manner in which $\beta Uw_x$ perturbs the in vacuo eigenmodes of the system; however, to view the flutter dynamics without seeing energies grow exponentially in time (as per the previous point in this section), one must include nonlinear restoring force $b_2>0$. 

\vskip.1cm
\noindent{\bf Instability in Each Critical Parameter}
\vskip.1cm
As can be expected, we can de/stabilize the dynamics through a myriad of parameters. Unsurprisingly, the flow parameters $\beta$ and $U$ are the principal candidates for the onset of instability. However, within certain limits, we can utilize $L$ and $b_1$ (for {\bf C} and {\bf H}) to affect the onset of instability. On the other hand, if we impose damping of the form $k_0$, we can re-stabilize eigenvalues that have become unstable due to any of the aforementioned parameters. This is to say that damping of the form $k_0w_t$ (for $k_0$ sufficiently large) is ``strong? enough to stabilize flutter. For a fixed $U > \Ucrit$, increasing $k_0$ provides two regimes: initially, it increases the rate of convergence to the stable LCO, but, for $k_0$ sufficiently large, the flutter dynamics are eventually damped out; see Figure \ref{fig11l}.

\vskip.1cm
\noindent{ \bf Relative Stability Across Configurations}
\vskip.1cm
As noted from various plots above, principally Figure \ref{fig1}, the configurations have a relative hierarchy of stability: the clamped panel {\bf C} is the most stable; the hinged panel {\bf H} is less stable; the cantilever {\bf CF} is least stable. This demonstrates the extent to which the boundary conditions (which serve to distinguish the in vacuo eigenfunctions) affect the inherent stability of each configuration. It is important to emphasize, then, that across all possible configurations one could consider with the boundary conditions presented here (including all left/right combinations), the {\bf CF} configuration is {\em the easiest} to make flutter.

\vskip.1cm
\noindent{\bf Uniformity of End Behavior}
\vskip.1cm
We have seen that, in most cases with no pre-stressing of the structure (when $b_1=0$), a very strong {\em uniformity} of end-behavior with respect to size and type (displacement, velocity) of initial data. This is true in the fully linear case, when trajectories are stable or when they are unstable, e.g., Figure \ref{fig4}; in the latter case, we note the same rate of exponential growth in time of trajectories independent of initial condition. In the nonlinear case---where interesting behaviors are found---we typically observe convergence to the same LCO for fixed parameters and varying initial conditions---see Figures \ref{fig11l} and \ref{fig23}.

\vskip.1cm
\noindent{\bf Reproduction of Qualitative Behaviors}
\vskip.1cm
It is worth emphasizing here that this {\em very} simplified model (extensible beam with piston theoretic RHS) for a very complex phenomenon (axial flow flutter) {\em can replicate all expected qualitative behaviors} from the engineering point of view. This includes convergence to simple LCOs (Section \ref{LCOs}), non-trivial LCOs (Figure \ref{fig12l}), and non-trivial (buckled) steady states (Figures \ref{fig8l} and \ref{fig9l}). We have even demonstrated a compelling case for deterministic chaos in Section \ref{chaosbaby}, by playing large flow values $U$ against large in-plane pre-stressing $b_1>>0$ for the hinged panel {\bf H}.

\appendix

\section{Fundamental Frequencies}

\begin{table}[H]
\begin{center}
\begin{tabular}{|l|r|r|}\hline
$n$& {\bf C};~~$k_nL$ & {\bf CF};~~ $k_nL$ \\\hline
1 &4.7300 & 1.8751 \\
2 &7.8532 & 4.6941\\
3 &10.9956 & 7.8548 \\
4 &14.1371 & 10.9955\\
5 &17.2787  & 14.1372\\
6 &20.4204 & 17.2788\\
7 &23.5619 & 20.4205\\
8 &26.7035 & 23.5619\\
9 & 29.8451 & 26.7035\\
10 & 32.9867&29.8451\\\hline
\end{tabular}
\end{center}
\caption{First 10 fundamental frequencies for the Clamped-Clamped ({\bf C})  and Clamped-Free ({\bf CF}) configuration.}
\label{table1}
\end{table}


\begin{thebibliography}{99}

\bibitem{pist2} Ashley, H. and Zartarian, G, 1956. Piston theory: A new aerodynamic tool for the aeroelastician. {\em Journal of the Aeronautical Sciences}, 23(12), pp. 1109--1118.

\bibitem{ball} Ball, J.M., 1973. Initial-boundary value problems for an extensible beam. {\em Journal of Mathematical Analysis and Applications}, 42(1), pp.61-90.

\bibitem{ball2} Ball, J.M., 1973. Stability theory for an extensible beam. {\em Journal of Differential Equations}, 14(3), pp.399-418.

\bibitem{slenderness} Han, S.M., Benaroya, H. and Wei, T., 1999. Dynamics of transversely vibrating beams using four engineering theories. {\em Journal of Sound and vibration}, 225(5), pp.935--988.

\bibitem{biacrip} Angela Cassia, Biazutti. and Helvecio Rubens, Crippa, 1994. Global attractor and inertial set for the beam equation. Applicable Analysis, 55(1-2), pp.61-78.

\bibitem{daniellorena} Bociu, L. and Toundykov, D., 2012. Attractors for non-dissipative irrotational von Karman plates with boundary damping. {\em Journal of Differential Equations}, 253(12), pp. 3568--3609.

\bibitem{bolotin}
Bolotin, V.V., 1963. {\em Nonconservative problems of the theory of elastic stability}. Macmillan.

\bibitem{che-tri:89:PJM} Chen, S.P. and Triggiani, R., 1989. Proof of extensions of two conjectures on structural damping for elastic systems. Pacific Journal of Mathematics, 136(1), pp.15-55.

\bibitem{quasi}
Chueshov, I., 2015. {\em Dynamics of Quasi-Stable Dissipative Systems}. Springer.


\bibitem{survey2} Chueshov, I., Dowell, E.H., Lasiecka, I. and Webster, J.T., 2016. Nonlinear Elastic Plate in a Flow of Gas: Recent Results and Conjectures. Applied Mathematics \& Optimization, 73(3), pp.475-500.

\bibitem{Memoires} Chueshov, I. and Lasiecka, I., 2008. {\em Long-time behavior of second order evolution equations with nonlinear damping}. American Mathematical Soc.

\bibitem{springer} Chueshov, I. and Lasiecka, I., 2010. {\em Von Karman Evolution Equations: Well-posedness and Long Time Dynamics}. Springer Science \& Business Media.

\bibitem{delay} Chueshov, I., Lasiecka, I. and Webster, J.T., 2014. Attractors for Delayed, Nonrotational von Karman Plates with Applications to Flow-Structure Interactions Without any Damping. Communications in Partial Differential Equations, 39(11), pp.1965--1997.

\bibitem{fereisl} Chueshov, I., Lasiecka, I. and Webster, J., 2014. Flow-plate interactions: Well-posedness and long-time behavior. Discrete \& Continuous Dynamical Systems-S, 7(5), pp.925-965.


\bibitem{dickey} Dickey, R.W., 1970. Free vibrations and dynamic buckling of the extensible beam. Journal of Mathematical Analysis and Applications, 29(2), pp.443-454.


\bibitem{dowellnon}
Dowell, E.H., 1966. Nonlinear oscillations of a fluttering plate I. {\em AIAA Journal}, 4(7), pp. 1267--1275. 1967. Nonlinear oscillations of a fluttering plate. II. {\em AIAA Journal}, 5(10), pp. 1856--1862.

\bibitem{dowellchaos} Dowell, E.H., 1982. Flutter of a buckled plate as an example of chaotic motion of a deterministic autonomous system. {\em Journal of Sound and Vibration}, 85(3), pp.333-344.


\bibitem{dowellapprox1} Jaworski, J.W. and Dowell, E.H., 2008. Free vibration of a cantilevered beam with multiple steps: Comparison of several theoretical methods with experiment. {\em Journal of Sound and Vibration}, 312(4-5), pp.713--725.

\bibitem{dowell} Dowell, E.H., Clark, R. and Cox, D., 2004. {\em A modern course in aeroelasticity} (Vol. 3). Dordrecht: Kluwer academic publishers.

\bibitem{inext1} Dowell, E. and McHugh, K., 2016. Equations of Motion for an Inextensible Beam Undergoing Large Deflections. Journal of Applied Mechanics, 83(5), p.051007.

\bibitem{edenmil} Eden, A. and Milani, A.J., 1993. Exponential attractors for extensible beam equations. Nonlinearity, 6(3), p.457.



\bibitem{han-las:00} Hansen, S.W. and Lasiecka, I., 2000. Analyticity, hyperbolicity and uniform stability of semigroups arising in models of composite beams. Mathematical Models and Methods in Applied Sciences, 10(04), pp.555-580.

\bibitem{holmes} P.J. Holmes, Bifurcations to divergence and flutter in flow-induced oscillations: a finite dimensional analysis. {\em Journal of Sound and Vibration}, 53(4), 1977, pp.471--503.

\bibitem{HolMar78} P. Holmes and J. Marsden, Bifurcation to divergence and flutter in flow-induced oscillations: an infinite dimensional analysis.
{\em Automatica}, 14, 1978,  pp.367--384.

\bibitem{HLW} Howell, J., Lasiecka, I., and Webster, J.T., 2016, Quasi-stability and Exponential Attractors for A Non-Gradient System---Applications to Piston-Theoretic Plates with Internal Damping, {\em EECT}. 5(4) pp. 567--603.

\bibitem{HTW} Howell, J.S., Toundykov, D. and Webster, J.T., 2018. A Cantilevered Extensible Beam in Axial Flow: Semigroup Well-posedness and Post-flutter Regimes. {\em SIAM Journal on Mathematical Analysis}, 50(2), pp.2048-2085.

\bibitem{modalcant} Huang, L. and Zhang, C., 2013. Modal analysis of cantilever plate flutter. {\em Journal of Fluids and Structures}, 38, pp.273--289.


\bibitem{inverted1} Kim, D., Coss\'{e}, J., Cerdeira, C.H. and Gharib, M., 2013. Flapping dynamics of an inverted flag. {\em Journal of Fluid Mechanics}, 736.

\bibitem{lagnese} Lagnese, J., 1989. {\em Boundary Stabilization of Thin Plates}. SIAM Studies in Applied Mathematics.

\bibitem{lagleug} Lagnese, J.E. and Leugering, G., 1991. Uniform stabilization of a nonlinear beam by nonlinear boundary feedback. {\em Journal of Differential Equations}, 91(2), pp. 355--388.


\bibitem{conequil1}
Lasiecka, I. and Webster, J., 2014. Eliminating flutter for clamped von Karman plates immersed in subsonic flows. {\em Communications on Pure \& Applied Analysis}, 13(5), pp. 1935--1969.
\newblock Updated version (May, 2015): {\url{http://arxiv.org/abs/1409.3308}}.

\bibitem{conequil2}
Lasiecka, I. and Webster, J.T., 2016. Feedback stabilization of a fluttering panel in an inviscid subsonic potential flow. {\em SIAM Journal on Mathematical Analysis}, 48(3), pp.1848-1891.


\bibitem{beam4} Ma, T.F., Narciso, V. and Pelicer, M.L., 2012. Long-time behavior of a model of extensible beams with nonlinear boundary dissipations. Journal of Mathematical Analysis and Applications, 396(2), pp.694-703.

\bibitem{b:paidoussis:98} Paidoussis, M.P., 1998. Fluid-structure interactions: slender structures and axial flow (Vol. 1). Academic press.


\bibitem{russell:93:JMAA} Russell, D.L., 1993. A general framework for the study of indirect damping mechanisms in elastic systems. Journal of Mathematical Analysis And Applications, 173(2), pp.339-358.

\bibitem{beamdamping} Russell, D.L., 1991. A comparison of certain elastic dissipation mechanisms via decoupling and projection techniques. Quarterly of Applied Mathematics, 49(2), pp.373-396.



\bibitem{paidoussis} Semler, C., Li, G.X. and Paidoussis, M.P., 1994. The non-linear equations of motion of pipes conveying fluid. {\em Journal of Sound and Vibration}, 169(5), pp.577-599.

\bibitem{amjad} Serry, M. and Tuffaha, A., 2017. Static Stability Analysis of a Thin Plate with a Fixed Trailing Edge in Axial Subsonic Flow: Possio Integral Equation Approach. arXiv preprint arXiv:1708.06956.

\bibitem{shubov1}
M. Shubov,  Solvability of reduced Possio integral equation in theoretical aeroelasticity, {\em Advances in Differential Equations}, 15, 2010,  pp.801--828.

\bibitem{inext2} Tang, D., Zhao, M. and Dowell, E.H., 2014. Inextensible beam and plate theory: computational analysis and comparison with experiment. Journal of Applied Mechanics, 81(6), p.061009.

\bibitem{jfs}
Vedeneev, V.V., 2013. Effect of damping on flutter of simply supported and clamped panels at low supersonic speeds. {\em Journal of Fluids and Structures}, 40, pp. 366--372.

\bibitem{vedeneev} Vedeneev, V.V., 2012. Panel flutter at low supersonic speeds. {\em Journal of fluids and structures}, 29, pp. 79--96.

\bibitem{webster} Webster, J.T., 2011. Weak and strong solutions of a nonlinear subsonic flow-structure interaction: Semigroup approach. {\em Nonlinear Analysis: Theory, Methods \& Applications}, 74(10), pp.3123-3136.

\bibitem{wonkrieg} Woinowsky-Krieger, S., 1950. The effect of an axial force on the vibration of hinged bars. {\em J. Appl. Mech.}, 17(1), pp.35--36.

\bibitem{paidoussis3} Zhao, W., Pa\"{i}doussis, M.P., Tang, L., Liu, M. and Jiang, J., 2012. Theoretical and experimental investigations of the dynamics of cantilevered flexible plates subjected to axial flow. Journal of Sound and Vibration, 331(3), pp.575-587.



\end{thebibliography}

\end{document}